\documentclass[review]{siamltex}
\usepackage{latexsym, graphicx, epsfig, amsmath, amsfonts,amssymb,subfigure}
\usepackage{multirow}
\usepackage{color}
 \usepackage[percent]{overpic}
\usepackage{algorithm}
\usepackage{algorithmicx}
\usepackage{algpseudocode}
\floatname{algorithm}{Algorithm}
\usepackage{threeparttable}
\usepackage{booktabs}

\def\mb{\mathbf}
\def\R{\mathbb{R}}
\def\V{\mathbb{V}}
\def\Nj{{N_{\text{obs}}}}
\def\Ns{{N_{\text{rv}}}}

\def\bigdot{\boldsymbol{\cdot}}
\def\calu{\mathcal{U}}
\def\bigdot{\boldsymbol{\cdot}}

\def\cale{\mathcal{E}}

\def\eref#1{{\rm (\ref{#1})}}

\title{A non-intrusive reduced basis EKI for  time-fractional diffusion inverse problems}%

\author{Fenglian Yang \thanks{College of Science, Hohai University, Nanjing 210098, China (yangfenglian@hhu.edu.cn). This author's work is partially supported by the NSF of China (No. 1160118).}
          \and Liang Yan \thanks{Department of Mathematics,  Southeast University, Nanjing, 210096, China (yanliang@seu.edu.cn). This author's work was supported by NSF of China (No.11771081), the science challenge project (No. TZ2018001),  Qing Lan project of Jiangsu Province and the Southeast University Zhishan Young Scholars Program. }}
\begin{document}
\graphicspath{figure/}
\maketitle

\begin{abstract}
In this study we consider an  ensemble Kalman inversion (EKI) for the numerical solution of time-fractional diffusion  inverse problems (TFDIPs).  Computational challenges in the EKI arise from the need for repeated evaluations of the forward model. We address this  challenge by introducing a non-intrusive reduced basis (RB) method for constructing surrogate models to reduce computational cost.  In this method, a reduced basis is extracted from a set of full-order snapshots by the proper orthogonal decomposition (POD), and a doubly stochastic  radial basis function (DSRBF) is used to learn the projection coefficients.  The DSRBF is carried out in the offline stage with a stochastic leave-one-out cross-validation algorithm to select the shape parameter, and the outputs for new parameter values can be obtained rapidly during the online stage.  Due to the complete decoupling of the offline and online stages, the proposed non-intrusive RB method -- referred to as POD-DSRBF -- provides a powerful tool to accelerate the EKI approach for TFDIPs. We demonstrate the practical performance of the proposed strategies through two nonlinear time-fractional diffusion inverse problems.  The numerical  results  indicate that the new algorithm can achieve significant computational gains without sacrificing accuracy.
\end{abstract}


\pagestyle{myheadings}
\thispagestyle{plain}
\markboth{Yang F. L. and Yan L.}
{ RB-EKI for TFDIPs }
\section{Introduction}
Consider the  following time-fractional  diffusion equation:
\begin{eqnarray}\label{fraceq}
\begin{cases}
^cD_t^{\alpha}u(x,t)=\nabla\cdot(\kappa(x) \nabla u(x,t))+f(x,t), &x\in \Omega, ~ t\in (0,T)\\
\mathcal{B}u(x,t)=g(x,t),&x \in \partial\Omega, t\in (0,T)\\
u(x,0)=u_0(x), &x\in \Omega,
\end{cases}
\end{eqnarray}
where  $\Omega$ is a bounded domain with  boundary $\partial\Omega$ and $T>0$ is a fixed value, $u$ represents the state variable,  $k(x)$ is the diffusion coefficient and  the right hand side $f$ denotes source terms. Here $\mathcal{B}$ is a boundary operator, $^c D^{\alpha}_t \,(0<\alpha<1)$ denotes the Caputo fractional derivative of order $\alpha$ with respect to $t$  and it is defined by 
\begin{eqnarray*}
^cD_t^{\alpha}u(x,t)=\frac{1}{\Gamma(1-\alpha)} \int ^t_0 \frac{\partial{u(x,s)}}{\partial{s}} \frac{ds}{(t-s)^{\alpha}}, \, 0<\alpha<1,
\end{eqnarray*}
where $\Gamma (\cdot)$ is the Gamma function. The system (\ref{fraceq}) has been widely used to model  anomalous subdiffusion in which the mean square variance grows slower than that in a Gaussian process, and has found a number of important practical applications \cite{Podlubny1999}.  The forward problems for time-fractional  diffusion equation (TFDE) have been extensively studied, see \cite{SY2011,Lin+Xu2007,FCY2013JCP} and references therein. However, in practice  there are many inputs unknown in the time fractional diffusion model, such as the fractional derivatives, the heat source and the diffusion coefficient and so on. These are time fractional diffusion inverse problems (TFDIPs). In general, the TFDIPs are much more difficult to solve  than the direct problem for TFDE due to the fact that the solution does not satisfy the general conditions of well-posedness \cite{JR2015}. There are, however, many works on both theoretical and numerical aspects over the last few years. We refer to the papers \cite{CNYY2009,JPY2017,LZJY2013,LYY2015,YY2015,ZX2011,ZW2010,TS2017,ZJY2018IP,LYY2015ANM,YY2014CMA} for a rather incomplete list, which range from recovering structural parameters of TFDEs to source terms and initial states.

In this paper, we consider the ensemble Kalman-based inversion method (EKI) \cite{Iglesias+Law+Kody2013ensemble,Iglesias2016regularizing,Law+Stuart2015data} to solve TFDIPs. The EKI  is a derivate-free approach that lies at the interface between the optimization and Bayesian approaches \cite{SS2017SJNA,SS2018AA}.  Since it only requires the evaluation of the forward operator but not its derivative, this approach has been successfully used to solve PDE-constrained parameter identification problems \cite{Iglesias2016regularizing}.  A potential drawback in the application of EKI to inverse problems is it requires a sufficiently large ensemble size to guarantee reliable estimations.  This is undesirable to solve the TFDIPs as each ensemble realization requires a solution of the forward model of TFDE. It is well known that the simulation of fractional-order systems is computationally demanding due to their long-range history dependence.  To address this challenge, surrogate modeling has been extensively developed for decades,  see \cite{Asher2015review,Li+Lin+Zhang2014adaptive,Li+Zhang+Lin2015surrogate,ju2018adaptive,YZ2018IJUQ}.  The key idea is to replace the high-fidelity full order system with a carefully constructed  computationally cheap approximated model, to reduce the memory needs and CPU time. Specifically,  the surrogate  model is often constructed offline and subsequently used online when running the EKI algorithm.  

Projection-type reduced order models (ROM)  are one possible realization of this idea. Here, the governing equations of the forward model are projected onto a limited number of reduced basis (RB) functions, which can be calculated via the proper orthogonal decomposition \cite{QMN2015book} or with greedy algorithm \cite{RHP2007reduced}, to obtain an inexpensive reduced order model.  During the online stage, for a new parameter value, the reduced model is constructed as a linear combination of the pre-computed RB functions, where the expansion coefficients are computed by projecting the full-order equation onto the RB space.  However, such scheme is problem-dependent and of an intrusive nature. Additionally, it causes computational inefficiency for complex nonlinear problems. To tackle these  challenges, non-intrusive methods \cite{HU2018JCP,GH2018CMAME} have been developed to construct a surrogate of the reduced coefficient so that the reduced coefficient can be recovered without requiring a projection of the full-order model. In this method, the data-fit surrogate models are constructed using interpolation and regression techniques to directly learn the map between the input parameters and the POD expansion coefficients.  There is also substantial past work on ROM for time-dependent problems, e.g. the non-intrusive frameworks based on radial basis function interpolation \cite{ADN2013NMPDE,XFPN2017CMAME}, and the Gaussian process regression \cite{GH2019CMAME}.  One can refer to \cite{BGW2015SIRV} for a comprehensive survey. 

Motivated by the recent developments in reduced order modeling, we propose a new non-intrusive RB method, namely POD-DSRBF, for constructing surrogate models to accelerate the EKI for TFDIPs.  After extracting the RB functions from a set of snapshots by POD, a learning approach using a doubly stochastic radial basis function (DSRBF) is used to establish a mapping from parameter values to projection coefficients onto the RB space.  Equipped with a stochastic leave-one-out cross-validation (LOOCV) parameter selection for the shape parameter of the RBF, the efficiency of the DSRBF can be further enhanced.   We shall discuss the basic idea and the efficient implementation of the algorithms.
To the best of our knowledge, this is the first investigation of  the EKI combining with the non-intrusive RB method to solve TFDIPs.  Numerical results indicate that the RB-EKI approach shows good performance in both accuracy and efficiency for TFDIPs.

The structure of the paper is as follows. In the next section, we review the formulation of  time fractional diffusion 
inverse problems and the solution approach via EKI.  In section 3, we introduce the POD-DSRBF approach to construct accurate surrogate models. In section 4, we use two different time fractional diffusion inverse problems to demonstrate the accuracy and efficiency of our method. Finally, we conclude the paper with Section 5. 

\section {Background and problem formulation}\label{sec:setup}
In this section, we first give a brief overview of the time fractional diffusion inverse problems (TFDIPs). Then we will introduce the EKI  to solve TFDIPs. 
\subsection{Mathematical formulation}
Unlike integer order differential equations, the analytical solutions of time fractional diffusion equation are usually not available.  We consider a numerical discretization of the system (\ref{fraceq}), described by 
\begin{equation}\label{stateq}
A(u_h(t; \theta), \theta) = 0,
\end{equation}
where $u_h(t;\theta): \mathcal{T} \times \Theta \rightarrow  \R^{n_h}$ is the discrete solution, $n_h$ is the dimension of the finite-dimensional discretization of the spatial domain, and $\theta \in \Theta \subset \R^{d}$ is the  $d$-dimensional parameter vector. The discrete operator $A$ denotes a numerical approximation, e.g., by the finite element  or finite difference method. Typically, the $n_h$ is large in order to resolve the details of the system with high accuracy.  The goal of TFDIPs is to estimate the unknown parameters $\theta$ from noisy observations of the states $u_h$ given by
\begin{equation}\label{dataeq}
y_{obs} = g(u_h; \theta)+\xi.
\end{equation}
Here $g$ is a discretized observation operator mapping from the states and parameters to the observables, $\xi \in \R^m$ is the measurement error.  We assume that the error $\xi$ is a Gaussian random vector with mean zero and covariance matrix $\Gamma \in \R^{m\times m}$, i.e., $\xi \sim N(0, \Gamma)$.
 The system model (\ref{stateq}) and observation model (\ref{dataeq}) together define a forward model $y=f(\theta)$ that maps the unknown parameter to the observable data. 
 
\subsection{Ensemble Kalman inversion}
Ensemble Kalman inversion (EKI) is a recently proposed inversion methodology that lies at the interface between the optimization and Bayesian approaches \cite{Iglesias+Law+Kody2013ensemble,SS2017SJNA}. It applies the ensemble Kalman filter (EnKF) \cite{Evensen1994} to the inverse problem setting by introducing a trivial dynamics for the unknown parameters.  In this work, we follow closely  the framework in \cite{Iglesias2016regularizing}. 

To formulate the TFDIPs in a Bayesian framework, we model the parameter $\theta$ as a random variable, endow it with a prior distribution $\pi(\theta)$.   Assume we derive $N_e$ initial ensemble $\theta_0^{(j)}\, (j \in \{1, \cdots, N_e\})$ from the prior $\pi(\theta)$.  The algorithm works by iteratively updating an ensemble of candidate solutions $\{\theta^{(j)}_{n}\}^{N_e}_{j=1}$ from iteration index $n$ to $n+1$.   Define the ensemble mean
$$\bar{\theta}_n=\frac{1}{N_e}\sum^{N_e}_{j=1}\theta^{(j)}_n,\, \bar{\omega}_n =\frac{1}{N_e}\sum^{N_e}_{j=1}f(\theta^{(j)}_n),$$
and convariances
\begin{equation}\label{coveq1}
C^{\theta \omega}_{n} = \frac{1}{N_e-1}\sum^{N_e}_{j=1}(\theta^{(j)}_n-\bar{\theta}_n)(f(\theta^{(j)}_n)-\bar{\omega}_n)^T, 
\end{equation}

\begin{equation}\label{coveq2}
C^{\omega \omega}_{n} = \frac{1}{N_e-1}\sum^{N_e}_{j=1}(f(\theta^{(j)}_n)-\bar{\omega}_n)(f(\theta^{(j)}_n)-\bar{\omega}_n)^T.  
\end{equation}

Then the EKI update formulae are given by
\begin{equation}\label{upEn}
\theta^{(j)}_{n+1}=\theta^{(j)}_{n}+C^{\theta\omega}_n(C_n^{\omega\omega}+\gamma_n \Gamma)^{-1}(y^{(j)}_{n+1}-\omega^{(j)}_n),
\end{equation}
where 
\begin{equation*}
y^{(j)}_{n+1} = y_{obs}+ \xi^{(j)}_{n+1}, \,\,   \xi^{(j)}_{n+1}\sim N(0, \Gamma),
\end{equation*}
and the regularization parameter $\gamma_n$ is choosing by 
\begin{equation*}\label{alrule}
\gamma^N_n \|\Gamma^{1/2}(C_n^{\omega\omega}+\gamma^N_n\Gamma)^{-1}(y^{(j)}-\bar{\omega}_n)\|\geq \rho \|\Gamma^{-1/2}(y^{(j)}-\bar{\omega}_n)\|.
\end{equation*}

As an iterative regularization method, the  iterative EKI is terminated according to the following discrepancy principle
\begin{equation}\label{dprule}
\|\Gamma^{-1/2}(y_{obs}-\bar{\omega}_n)\| \leq \tau  \|\Gamma^{-1/2}(y_{obs}-f(\theta^{\dag}))\|,
\end{equation}
where $\tau$ is a constant and  $\theta^{\dag}$ denotes the truth properties.  

Noted that the  total cost of an $N_e$ size ensemble of EKI is approximately $N_eJ$ forward model evaluations where $J$ is the total number of iterations. Therefore, each iteration could involve over $N_e$ TFDE solvers even for relatively small $J$. Hence,  the cost of EKI for TFDIPs is prohibited as it requires a relatively large ensemble size to guarantee its performance.  It is thus natural to construct a surrogate of the forward model before the data are available. In the next section, we will focus on the reduced basis method, which is widely used in applied mathematics and engineering.

\section{The reduced basis EKI method}\label{sec:method}

\subsection{The proper orthogonal decomposition and the reduced basis space}

The goal of model order reduction is to reduce the high dimensional model (\ref{stateq}), in the sense of reducing the number of degrees of freedom, with the expectation of computational savings in terms of both storage and CPU times.  One of the convenient tools for model order reduction is the reduced basis (RB) method. It seeks the approximate solution to (\ref{stateq}) in a reduced space spanned by a set of parameter-independent RB functions.

To generate an RB space, we consider the set of $Q=n_t n_{\theta}$ snapshots, $\{u_h(t_i;\theta_j)| i=1,\cdots, n_t, j=1,\cdots, n_{\theta} \}$, which are snapshots at $n_t$ different time instances $t_1,\cdots, t_{n_t} \in \mathcal{T} $ and $n_{\theta}$ different inputs  $\theta_1,  \cdots, \theta_{n_{\theta}} \subset \Theta$.  Define the snapshot matrix $\mb{S} \in \R^{n_h\times Q}$, which  collects the snapshots $u_h(t_i;\theta_j)$ as its columns. Thus, each row in the snapshot matrix corresponds to a spatial location and each column corresponds to a snapshot.

 With this snapshot matrix, a (thin) singular value decomposition (SVD) is then performed to obtain the reduced space:
 \begin{equation*}
\mb{S} =\mb{U\Sigma V}^T,
\end{equation*}
 where $\mb{U}\in \R^{n_h\times Q}$ and $\mb{V}\in \R^{Q\times Q}$ are the left and right singular vectors of $\mb{S}$ respectively. The singular values $\sigma_1 \geq \sigma_2\geq \cdots \geq \sigma_Q \geq 0$ of $\mb{S}$ give the diagonal matrix $\mb{\Sigma}=\mbox{diag}\{ \sigma_1,\sigma_2, \cdots,\sigma_Q\} \in \R^{Q\times Q}$.
The POD basis  $\mb{U}_p = [\mb{u}_1, \mb{u}_2,\cdots, \mb{u}_p]\in \R^{n_h\times p}$ is  defined as the $p$ left singular vectors of $\mb{S}$ that correspond to the $p$ largest singular values. Then the RB space $\V_p$ is defined as 
\begin{equation*}
\V_p =\mbox{Span} \{\mb{u}_1, \mb{u}_2,\cdots, \mb{u}_p\}.
\end{equation*}
Among all orthonormal bases of size $p$, the POD basis minimizes the  error of snapshot reconstruction
\begin{equation*}
\min_{\mb{U}_p \in \R^{n_h\times p}} \|\mb{S}-\mb{U}_p\mb{U}_p^T\mb{S}\|^2_F =\sum^{Q}_{k=p+1}\sigma_k^2,
\end{equation*}
where $\|\cdot\|_F$ denotes the Frobenius norm.
Thus, the singular values provide a guidance for choosing the size of the POD basis. 
A typical approach is to choose $p$ such that 
\begin{equation*}
\frac{\sum^p_{i=1}\sigma_i^2}{\sum^Q_{i=1}\sigma_i^2}>\epsilon_{pod},
\end{equation*}
where $\epsilon_{pod}$ is a user-specified tolerance. 

Once the RB space is available, the field $u_h$ can be approximated by a linear expansion in the POD basis:
\begin{equation*}
u_p(t; \theta) = \sum^p_{k=1} \mb{u}_k a_k(t;\theta),
\end{equation*}
where $a_k(t; \theta)$ denote the POD expansion coefficients.  One can evaluate the POD coefficient for a given parameter value $\theta$ by interpolation \cite{BJWW2000IP} or the Galerkin procedure \cite{QMN2015book}. In this study, we propose a learning approach, namely doubly stochastic radial basis function (DSRBF), to calculate the POD expansion coefficients.   Note that, given a snapshot $u_h(t;\theta)$, we can compute its representation in the POD basis via the coefficients $a_k(t;\theta) = \mb{u}_k^T u_h(t;\theta), \, k=1,\cdots, p$, where we have used that the POD basis vectors are orthonormal. Our learning task is now transformed into learning a model for the POD coefficients $a_k(t;\theta)$ from a training data 

\begin{eqnarray}\label{tr_data}
\begin{array}{lr}
\mathcal{D} =\Big\{\{(t, \theta), \mb{U}^T u_h(t;\theta)\}: t\in \mathcal{T}_{tr}, \theta\in \Theta_{tr}\Big\},\\
\mathcal{T}_{tr}=\{t_1,\cdots,t_{{N^{tr}_{t}}}\}\subset \mathcal{T},\\
\Gamma_{tr}=\{\theta_1,\cdots,\theta_{{N^{tr}_{\theta}}}\}\subset \Theta.\\
 \end{array}
\end{eqnarray}

\subsection{Doubly stochastic radial basis function}

In this subsection, we briefly describe the doubly stochastic radial basis function (DSRBF) method.  We refer interested readers to \cite{YYL2018JCP} for more details of the development of the DSRBF strategy. 

In DSRBF,  we first pick a \emph{translation-invariant} radial kernel $K=\varphi(\|\bigdot-\bigdot\|_{\ell^2}):\R^d\times\R^d\to\R$ with a function $\varphi:\R\to\R$  known as the radial basis function (RBF).  Commonly used RBF, include the Gaussian $\varphi(r)=\exp(-r^2)$, the multiquadrics (MQ) $\varphi(r)=(1+r^2)^{\beta/2}$ with $\beta=1$  and inverse multiquadrics (IMQ) with $\beta=-1$. Then, we can define  the finite-dimensional \emph{doubly stochastic trial space}  in the form of
\begin{equation}\label{ds trial space}
  \calu_{Z,\cale}:= \text{Span}\{\varphi(\varepsilon_j\|\bigdot-z_j\|_{\ell^2}\;|\;z_j\in Z,\,\varepsilon_j\sim\cale\}
\end{equation}
for some quasi-uniform trial centers $Z$ and some stochastic shape parameters following the probability distribution $\cale$.

For learning problems using DSRBF, one seeks an approximant   $ u \in \calu_{Z,\cale}$ of an unknown function $f$ from a collection of examples, called \emph{training data}. Recall that any trial function is a linear combination of the basis used in defining \eref{ds trial space} and is in the form of
\begin{equation}\label{num exp}
    u(x) =\sum_{z\in Z}c_j \varphi_j(x) := \sum_{z\in Z}c_j \varphi(\varepsilon_j\|x-z_j\|_{\ell^2})
\end{equation}
with $\varepsilon_j\sim \cale$ and
for some coefficients $\mb{c}=[c_1,\ldots,c_{n_Z}]^T\in\R^{n_Z}$. To have a well-posed fully-discretized problem, we make an observation  $\varepsilon =\{\varepsilon_j\}_{j=1}^{n_Z}\sim \cale^{n_Z}$ and define the least-squares numerical solution
\begin{equation*}
  u_{X,Z,\varepsilon} := \arg\inf_{u\in\calu_{{Z,\cale}}} \sum_{x\in X} | u - f |^2.
\end{equation*}
Using the training data $\{(z_i, y_i): z_i \in Z, y_i=f(z_i), i=1,\ldots, n_Z \}$, we can  yield the square unsymmetric system
\begin{equation}\label{Aa=f}
    A_{\varepsilon}(Z,Z) \mb{c} = f(Z),
\end{equation}
where $\mb{c}$ is the vector of unknown coefficients, and $f(Z)$ is the data vector.  The matrix $A_{\varepsilon}(Z,Z)$ is an $n_Z\times n_Z$ coefficient matrix with entries $[A_{\varepsilon}(Z,Z)]_{ij} = \varphi(\varepsilon_j\|z_i-z_j\|_{\ell^2})$ for $z_i,z_j\in Z$.  Therefore, all unknown coefficients $\{c_j\}$ can be easily obtained by a standard matrix solver. Once $\{c_j\}$ have been determined, the approximation solution $u_{X,Z,\varepsilon}$ can be evaluated from \eref{num exp}.

\begin{algorithm}[t]
  \caption{Stochastic LOOCV cost vector}
  \label{alg:SLOOCV}
  \begin{algorithmic}[1]
    \Require
  Given the shape parameter  $\varepsilon$, the RBF $\varphi$ and the training set $\mathcal{D}=\{Z,f_Z\}$,  a small integer $\Ns<N_Z$.
 \State  Generate a sequence $\{\omega_j \in \R^{N_Z}\}_{j=1}^{\Ns}$ of $\Ns$ normal random vectors and compute $v_j = A_\varepsilon(Z,Z) w_j, \, j=1,\cdots, \Ns$.
 \State Compute $\beta_1= \Big[ \sum^\Ns_{k=1} v_k \odot w_k \Big] \oslash \Big[ \sum^\Ns_{k=1} v_k \odot v_k \Big]$, where $\odot $ and $\oslash$ represent element-wise multiplication and division operators of vectors, respectively.
\State Compute $\beta_2 = W V^+f_{Z}$, where $W=[w_1,\ldots,w_{\Ns}]$, $V=[v_1,\ldots,v_{\Ns}]$ and  $V^+$ is the pseudo-inverse of $V$. 
   \State
    \Return $e(\varepsilon) =\beta_2\oslash \beta_1$.
  \end{algorithmic}
\end{algorithm}

For efficiency, one should choose an `optimal' shape parameter $\varepsilon_j\sim \cale$.  In this work,  we use a stochastic leave-one-out cross validation (LOOCV) proposed in \cite{YYL2018JCP} to select this parameter.  The stochastic LOOCV method basically consists of the following two steps:
\begin{itemize}
  \item[(1)] For the linear equation $ A_{\varepsilon}(Z,Z) \mb{c} = f(Z)$, whose matrix depends on a single shape parameter $\varepsilon>0$, we collect a set of observations $\widetilde{\varepsilon} = [\widetilde{\varepsilon}_1,\ldots,\widetilde{\varepsilon}_\Nj]^T$ of the `optimal' shape parameter. 
  \item[(2)] Generate a random vector $[\varepsilon_j]_{j=1}^{n_Z} \sim \cale^{n_Z} = [ \chi^2( \text{mean}(\widetilde{\varepsilon}) )]^{n_Z}$.
\end{itemize}
The key ingredient of stochastic LOOCV is how to select the `optimal' shape parameter $\widetilde{\varepsilon}_j$. To this end,  we define a stochastic LOOCV cost vector $e(\varepsilon)=[e_1,\ldots,e_{n_Z}]^T$ using Algorithm \ref{alg:SLOOCV}. Then the stochastic-LOOCV optimal shape parameter can be defined as
\begin{equation}\label{e^*}
  \varepsilon^* = \arg\min_{\varepsilon>0} \| e(\varepsilon) \|_{\ell^2}.
\end{equation}
For $j=1$ to some $\Nj>0$, we apply \eref{e^*} to estimate its optimal shape parameter by some minimization algorithm. Then, store result as $\widetilde{\varepsilon}_j$.

We now present the DSRBF strategy. The DSRBF method consists of the following steps which are outlined here:
\begin{itemize}
  \item Generate a random vector $[\varepsilon_j]_{j=1}^{n_Z} \sim \cale^{n_Z} = [ \chi^2( \text{mean}(\widetilde{\varepsilon}) )]^{n_Z}$ using stochastic LOOCV.
  \item Construct the matrix system in \eref{Aa=f}, where
  \[
    [A_\varepsilon(Z,Z)]_{ij} = \varphi( \varepsilon_j(z_i - z_j) )
    \mbox{ and }
   [f(Z)]_j= f(z_j)
  \]
  for $1\leq i, j \leq n_Z$.
  \item Seek for numerical solution in the form of \eref{num exp} by solving \eref{Aa=f} for the expansion coefficient $\mb{c} \in\R^{n_Z}$.
\end{itemize}

\subsection{RB-based EKI algorithm}

\begin{algorithm}[t]  
  \caption{The offline and online stages for the POD-DSRBF method}  
  \label{alg:POD-DSRBF}  
  \begin{algorithmic}[1]  
  \State  {\bf Offline stage:} 
    \State Compute $Q$ full-order snapshots $\{u_h(t_i; \theta_j)\}$ and form the snapshot matrix $\mb{S}\in \R^{n_h\times Q}$;
    \State  Perform POD for $\mb{S}$ and get the $p$ orthogonal bases $\mb{U}_p \in \R^{n_h\times p}$;
    \State Prepare the training set $\mathcal{D} = \Big\{\{(t, \theta), \mb{U}^T u_h(t;\theta)\}: t\in \mathcal{T}_{tr}, \theta\in \Theta_{tr}\Big\}$;
    \State Construct the DSRBF model $\Phi(t, \theta; \epsilon)$ from $\mathcal{D}$.
    \vspace{0.3cm}
    \State {\bf Online stage:}
    \State Recover output $\Phi(t, \theta^*;\epsilon)$ for a new parameter value $\theta^*$;
    \State Evaluate the reduced-solution $\tilde{u}_p(t, \theta^*) = \mb{U}_p \Phi(t,\theta^*;\epsilon)$.
  \end{algorithmic}  
\end{algorithm}  

In this subsection, we will use the DSRBF combining with a tensor decomposition \cite{GH2019CMAME} to  learn the POD expansion coefficients.   For the $k$th entry $a_k=\mb{u}_k^T u_h$ of the projection coefficients, the training data can be written in a matrix as
\begin{equation*}
\mb{Q}_k = [a_k(t_i;\theta_j)]_{i,j}, 1\leq i\leq N^{tr}_t,\, 1\leq j \leq N^{tr}_{\theta},
\end{equation*}
as a result of the tensor grid between the time and the parameter locations in the training data.

We employ the SVD again to decompose the data of an expansion coefficient into several time- and parameter-modes
\begin{equation*}
\mb{Q}_k \approx \tilde{\mb{Q}}_k =\sum^{q_k}_{l=1} \lambda^k_{l} \Psi_l^k(\Phi^k_l)^T.
\end{equation*}
Here $\Psi_l^k$ and $\Phi^k_l$ are the $l$th discrete time- and parameter-modes for the $k$th projection coefficient, respectively, $\lambda^k_{l}$ is the $l$-th singular value for the same coefficient, and $q_k\,(1\leq q_k \leq p)$ is the corresponding truncation number.

DSRBFs are now trained to approximate the corresponding continuous modes as
\begin{eqnarray*}\label{tr_model}
\begin{array}{lr}
t \mapsto \tilde{\Psi}_l^k(t), \, \text{ trained from} \{(t_i,\Psi_l^k(t_i) ): i=1,\cdots, N^{tr}_t\} ,\\
\theta \mapsto \tilde{\Phi}^k_l(\theta), \, \text{trained from} \{(\theta_i,\Phi_l^k(\theta_i) ): i=1,\cdots, N^{tr}_{\theta}\},\\
 \end{array}
\end{eqnarray*}
where $\tilde{\Psi}_l^k(t)$ and $\tilde{\Phi}^k_l(\theta)$ are the $l$th continuous time and parameter DSRBF modes for the $k$th expansion coefficient, respectively. Then, a continuous learning function $\tilde{a}_k(t,\theta)$ for the $k$th expansion coefficient $a_k$ can be recovered as
\begin{equation*}
a_k(t; \theta) \approx \tilde{a}_k(t; \theta) =\sum^{q_k}_{l=1} \lambda^k_{l} \tilde{\Psi}_l^k(\tilde{\Phi}^k_l)^T, \,\, (t,\theta)\in \mathcal{T}\times \Theta.
\end{equation*}

The resulted DSRBF approximation of the POD coefficient vector $a(t; \theta)$ can  be also represented as follows:
\begin{equation}\label{RBFtrain}
\tilde{a}(t,\theta) = \Phi(t, \theta; \varepsilon),
\end{equation}
where $\Phi(t,\theta; \varepsilon)$ is the trained DSRBF model for the POD coefficient $a(t; \theta)$ and $\varepsilon$ is the shape parameter of the DSRBF.  The complete POD-DSRBF algorithm  is displayed in Algorithm \ref{alg:POD-DSRBF}.

\begin{algorithm}[t]
  \caption{ RB-EKI algorithm for  inverse problems}
  \label{alg:RBEKI}
  \begin{algorithmic}[1]
   \State  \textbf{Prior ensemble and perturbed noise.} Let $\rho<1$ and $\tau\geq1/\rho$. Generate
$$\theta^{(j)}_0\sim \pi(\theta), y^{(j)}=y_{obs}+\xi^{(j)}, \quad \xi^{(j)}\sim N(0,\Gamma), j=1,\cdots, N_e.$$
\noindent Then  for $n=1,2 \dots,$
    \State \textbf{Prediction step:}  Evaluate
\begin{equation*}\label{fevaluate}
\omega^{(j)}_n=\tilde{f}_p(\theta^{(j)}_n),\quad j=1,\cdots, N_e
\end{equation*}
and define $\bar{\omega}_n=\frac{1}{N_e}\sum^{N_e}_{j=1}\omega^{(j)}_n.$
    \State \textbf{Discrepancy principle:} If the discrepancy principle (\ref{dprule}) is satisfied,
stop. Output $\bar{\theta}_n=\frac{1}{N_e}\sum^{N_e}_{j=1}\theta^{(j)}_n.$
  \State \textbf{Analysis step:} Define $C^{\theta \omega}_{n}, C^{\omega \omega}_{n}$ using (\ref{coveq1}) and (\ref{coveq2}), respectively. 
\noindent Update each ensemble member using (\ref{upEn}).
  \end{algorithmic}
\end{algorithm}

Once the DSRBF is trained, we can predict the reduced coefficient $\tilde{a}(t; \theta^{*})$ for a new given parameter $\theta^{*}$, during the online stage. The corresponding reduced solution is given by
\begin{equation}\label{statePOD}
\tilde{u}_p(t,\theta^{*}) = \mb{U}_p \tilde{a}(t,\theta^{*}),
\end{equation}
where $\mb{U}_p$ is the POD basis set. Then the associated model outputs are 
\begin{equation}\label{dataPOD}
y_p = g(\tilde{u}_p(t,\theta^{*});\theta^{*}).
\end{equation}
If $p\ll n_h$, the dimension of the unknown state in \eref{statePOD} and \eref{dataPOD} are greatly reduced compared with that of the original full model \eref{stateq} and \eref{dataeq}. Thus, \eref{statePOD} and \eref{dataPOD} define a reduced-order model $y_p = \tilde{f}_p(\theta)$ that maps the parameter $\theta$ to an approximation of the model outputs $y_p$. 

It is clear that after obtaining the  reduced-order model, we can then replace the forward model $f$ by its approximation $\tilde{f}_p$, and obtain the RB-based EKI algorithm.  The pseudo-code of the RB-based EKI algorithm is presented in Algorithm \ref{alg:RBEKI}.  Notice that the computational cost of generating $N_e$ samples using $\tilde{f}_p$ in prediction step requires nothing but sampling of the DSRBF model of (\ref{RBFtrain}) with $N_e$ samples of $\theta$. This cost is minimal because it does not require any simulations of the forward model. 

\section{Numerical Examples}\label{sec:tests}
In order to illustrate the accuracy and efficiency of the RB-EKI approach for solving the time fractional diffusion inverse problems, in this section, we present numerical experiments with two type different inverse problems.  The first example, adapted from \cite{YZ2019JCP,yan+guo2015}, considers a  heat source inversion problem.  The second example is the problem of inferring the spatially-varying diffusion coefficient \cite{YZ2018IJUQ,Cui2014data}.  For numerical examples, we set $\Omega=[0, 1]^2$ and $T=1$.   In all our tests, unless otherwise specified, we shall use the following parameters  $ \rho=1/\tau=0.7, \Ns =15, \Nj=10$.   We employ the MQ function in the DSRBF and the MATLAB function \texttt{fminbnd}  to solve \eref{e^*}. All the computations are performed using MATLAB 2015a on an Intel-i5 desktop computer.

\subsection{Example 1: 2D heat source inversion}
In this example, we consider the following model 
\begin{eqnarray}\label{2dsource}
\begin{array}{rl}
^cD_t^{\alpha}u-\nabla ^2 u&=e^{-t}\exp\Big[-0.5\Big(\frac{\|\theta-x\|}{0.1}\Big)^2\Big],\quad \Omega \times [0, 1],\\
\nabla u \cdot \textbf{n}&=0, \quad \mbox {on} \,\,\, \partial{\Omega},\\
u(x,0)&=0, \quad  \mbox{in}\,\,\, \Omega.
 \end{array}
\end{eqnarray}

The aim is to determine the source location $\theta=(\theta^1,\theta^2)$  from noisy measurements of the $u$-field at a finite set of locations and times.  The prior distributions on $\theta^i$ are independent and uniform, i.e., $\theta^i \sim U(0,1)$.  Below, unless stated otherwise, the observations are generated by selecting the values of the states at a uniform $3 \times 3$ sensor network.  At each sensor location, three measurements are taken at time $t = \{0.25, 0.75, 1\}$, which corresponds to a total of  27 measurements.  For any given values of $\theta$, the Eq. \eref{2dsource} is solved by 
using a finite difference/ spectral approximations (\cite{Lin+Xu2007}) with time step $\Delta t=0.01$ and polynomial degree $P=6$.  In order not to commit an `inverse crime',  we generate the data by solving the forward problem at a much higher resolution than that used in the inversion, i.e., with  $P=10$. In the examples below, unless otherwise specified, the Caputo fractional derivative of order $\alpha$ is 0.5.   

\subsubsection{Solution to the forward problem}

In this subsection, we assess the accuracy of the reduced order models for the forward model. We  construct a reduced order model using POD-DSRBF  whose solution captures the output of the forward model over the support of the prior distribution.  It should be noted that at this stage, measurement data do not enter the inference procedure. Therefore, we can construct the surrogate model off-line and save it for future use.

\begin{figure}
\begin{center}
  \begin{overpic}[width=0.5\textwidth,trim=20 0 20 15, clip=true,tics=10]{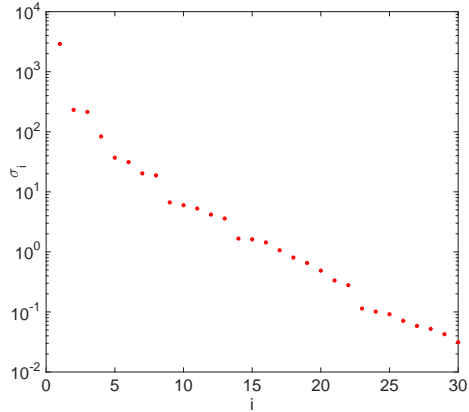}
    \end{overpic}
\end{center}
\caption{The first 30 singular values of the snapshot matrix $\mb{S}$.  }\label{sigvale}
\end{figure}

\begin{figure}
\begin{center}
  \begin{overpic}[width=0.45\textwidth,trim=20 0 20 15, clip=true,tics=10]{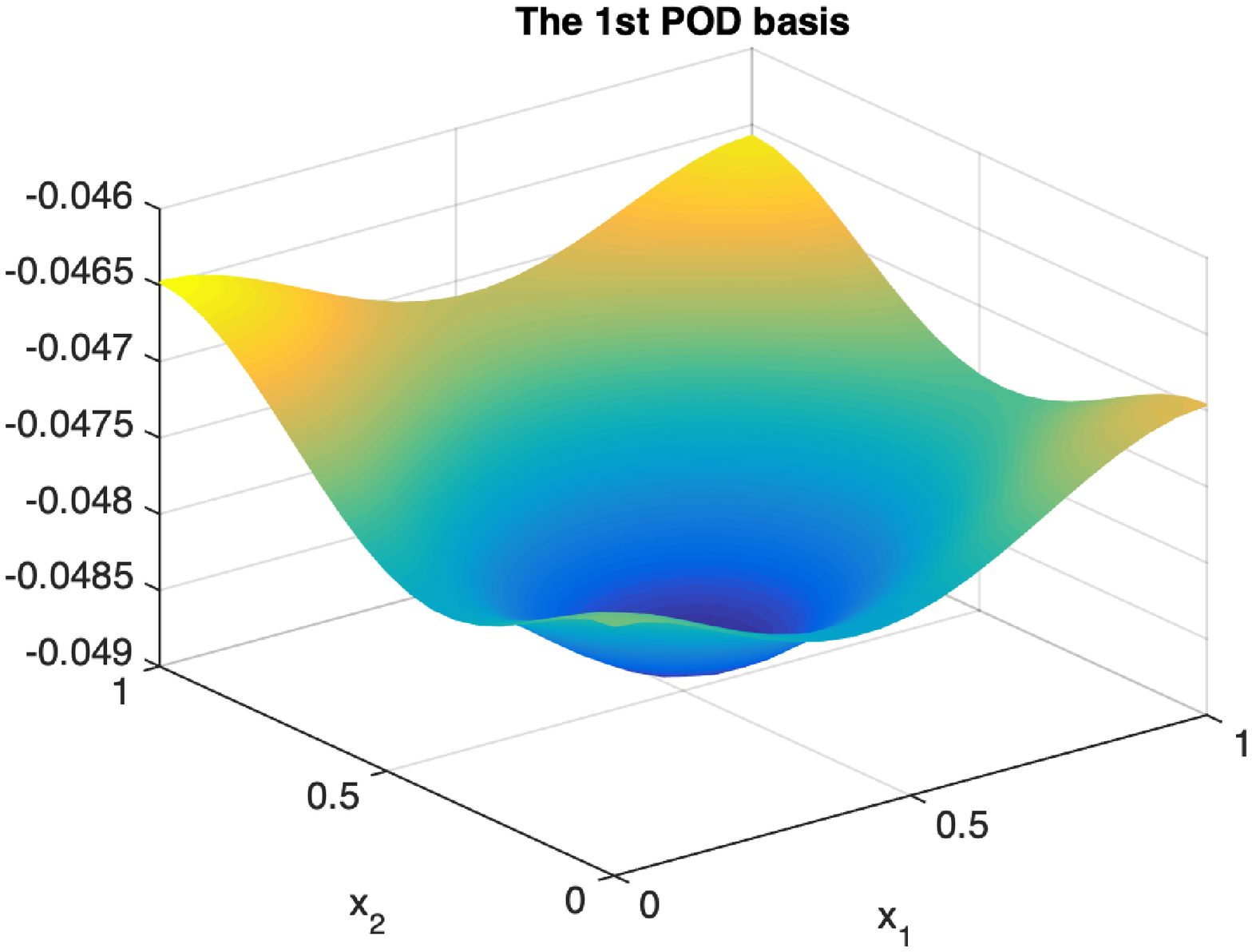}
  \end{overpic}
    \begin{overpic}[width=.45\textwidth,trim= 20 0 20 15, clip=true,tics=10]{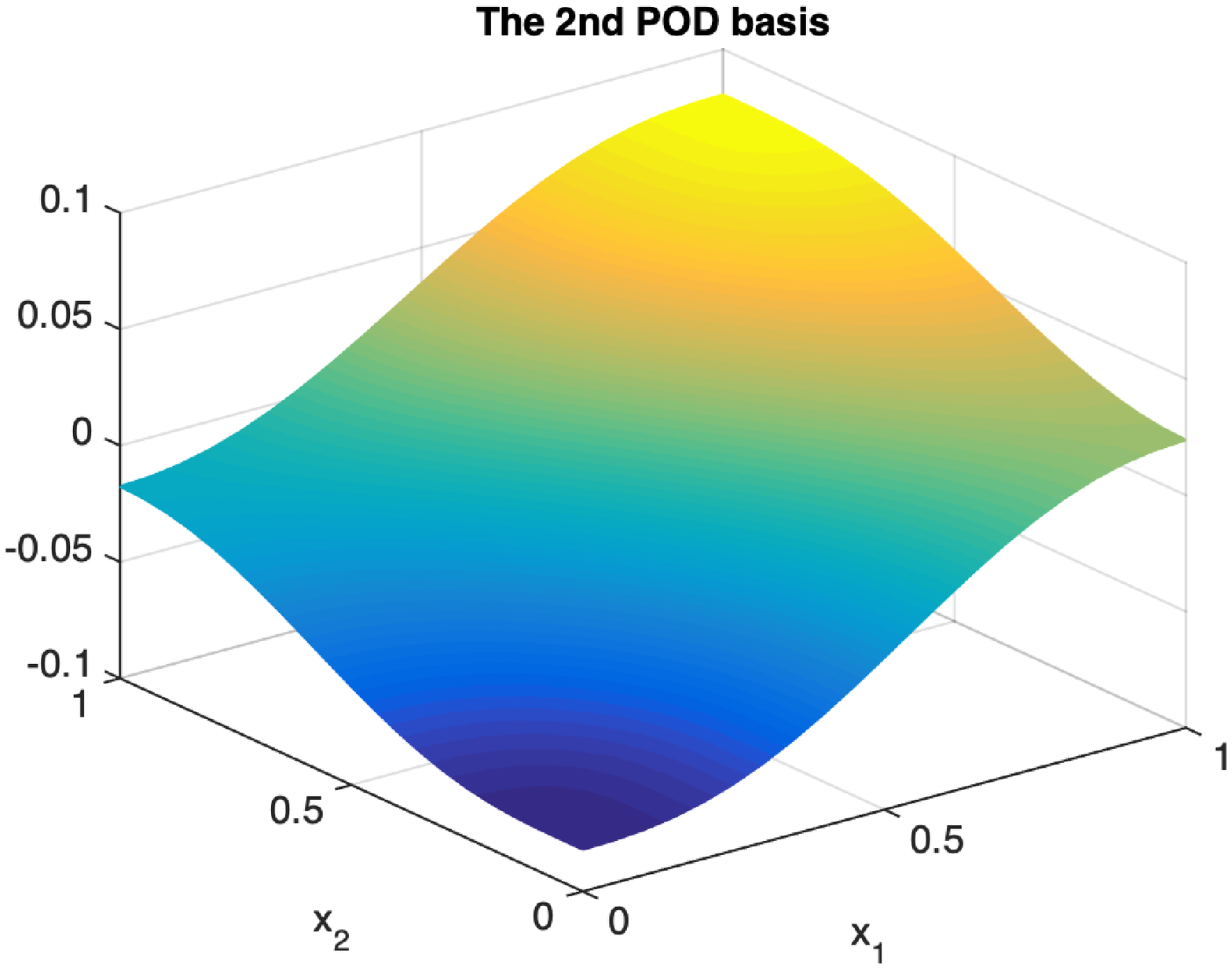}
  \end{overpic}
    \begin{overpic}[width=0.45\textwidth,trim= 20 0 20 15, clip=true,tics=10]{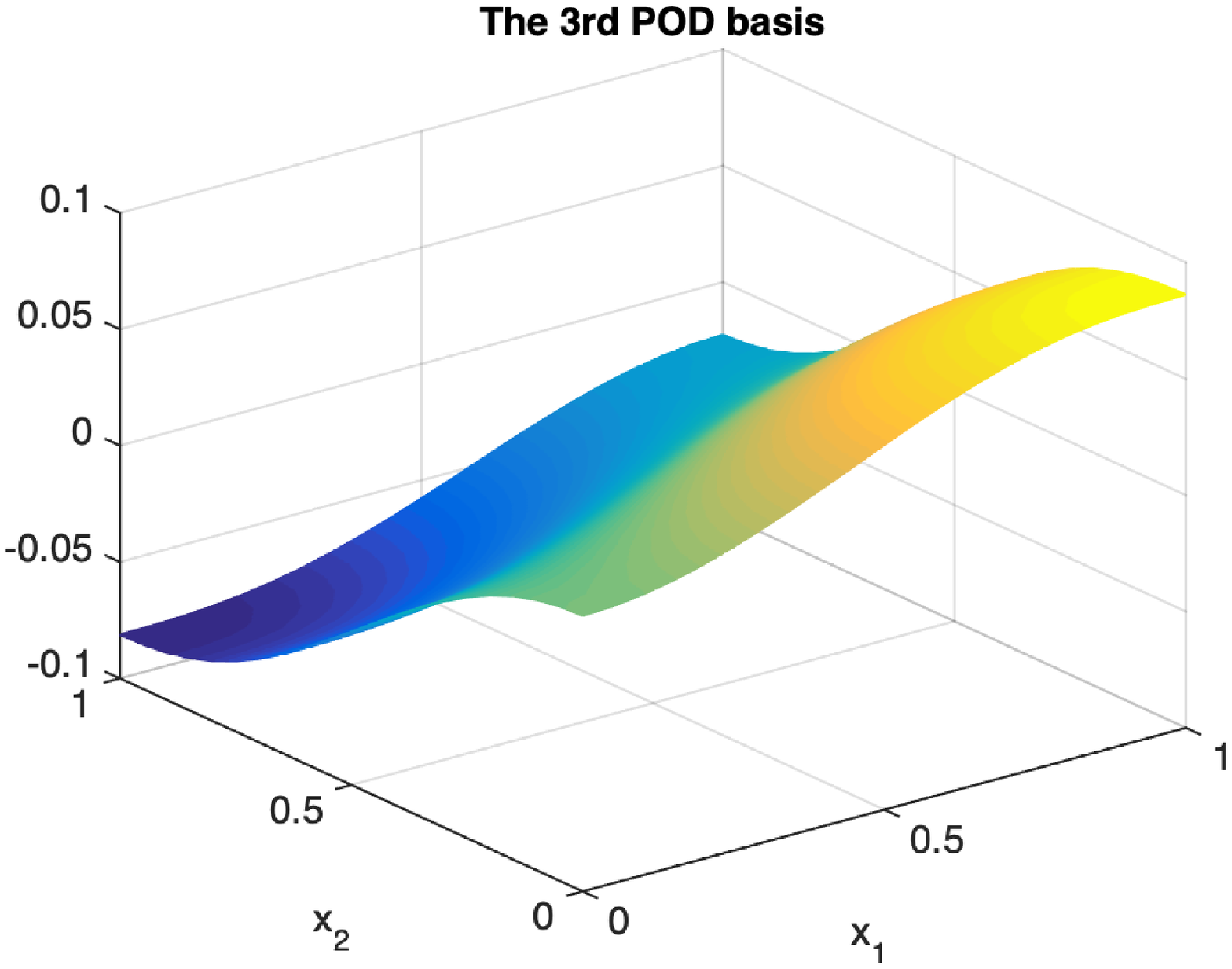}
  \end{overpic}
    \begin{overpic}[width=0.45\textwidth,trim=20 0 20 15, clip=true,tics=10]{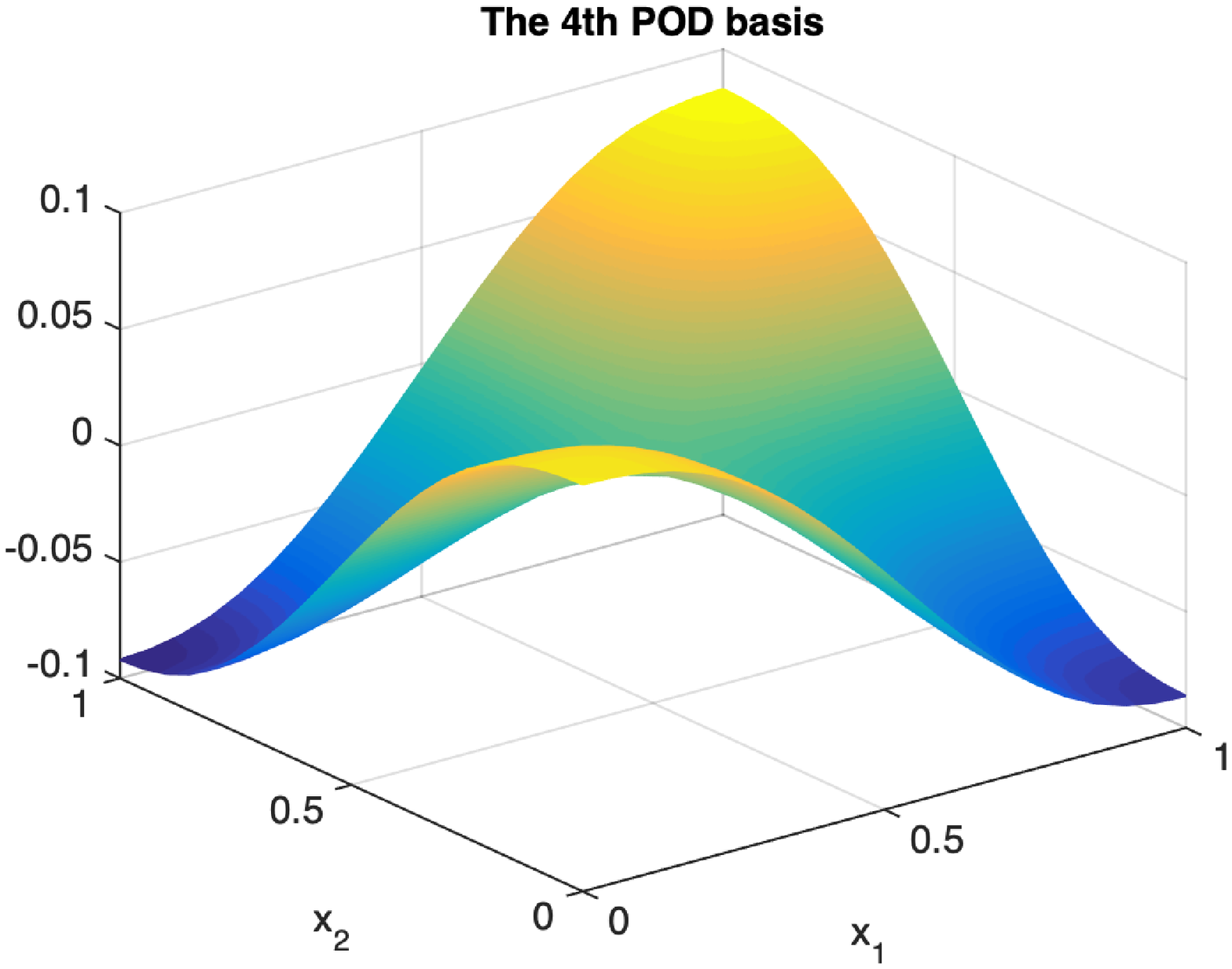}
  \end{overpic}
    \begin{overpic}[width=0.45\textwidth,trim= 20 0 20 15, clip=true,tics=10]{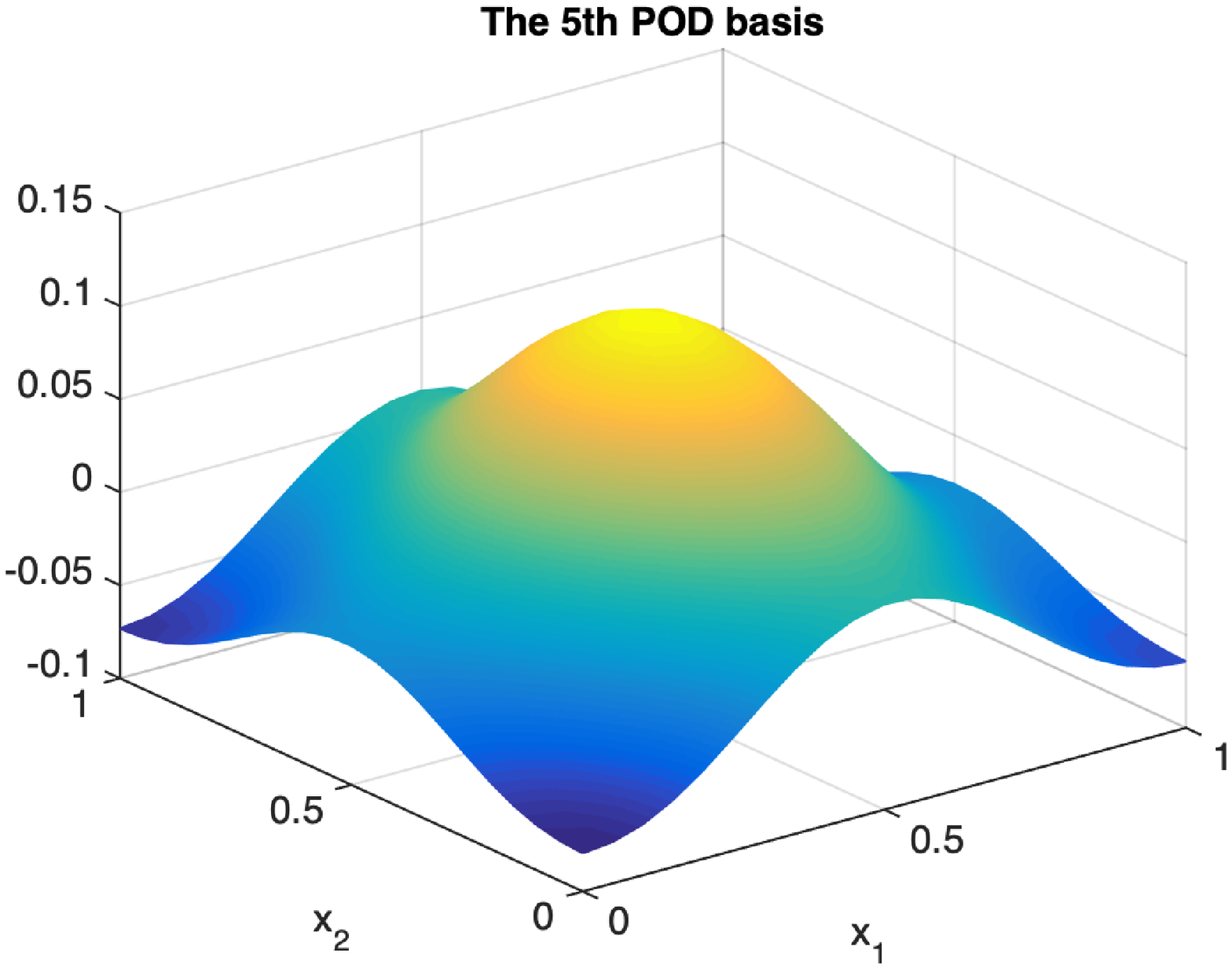}
  \end{overpic}
    \begin{overpic}[width=0.45\textwidth,trim= 20 0 20 15, clip=true,tics=10]{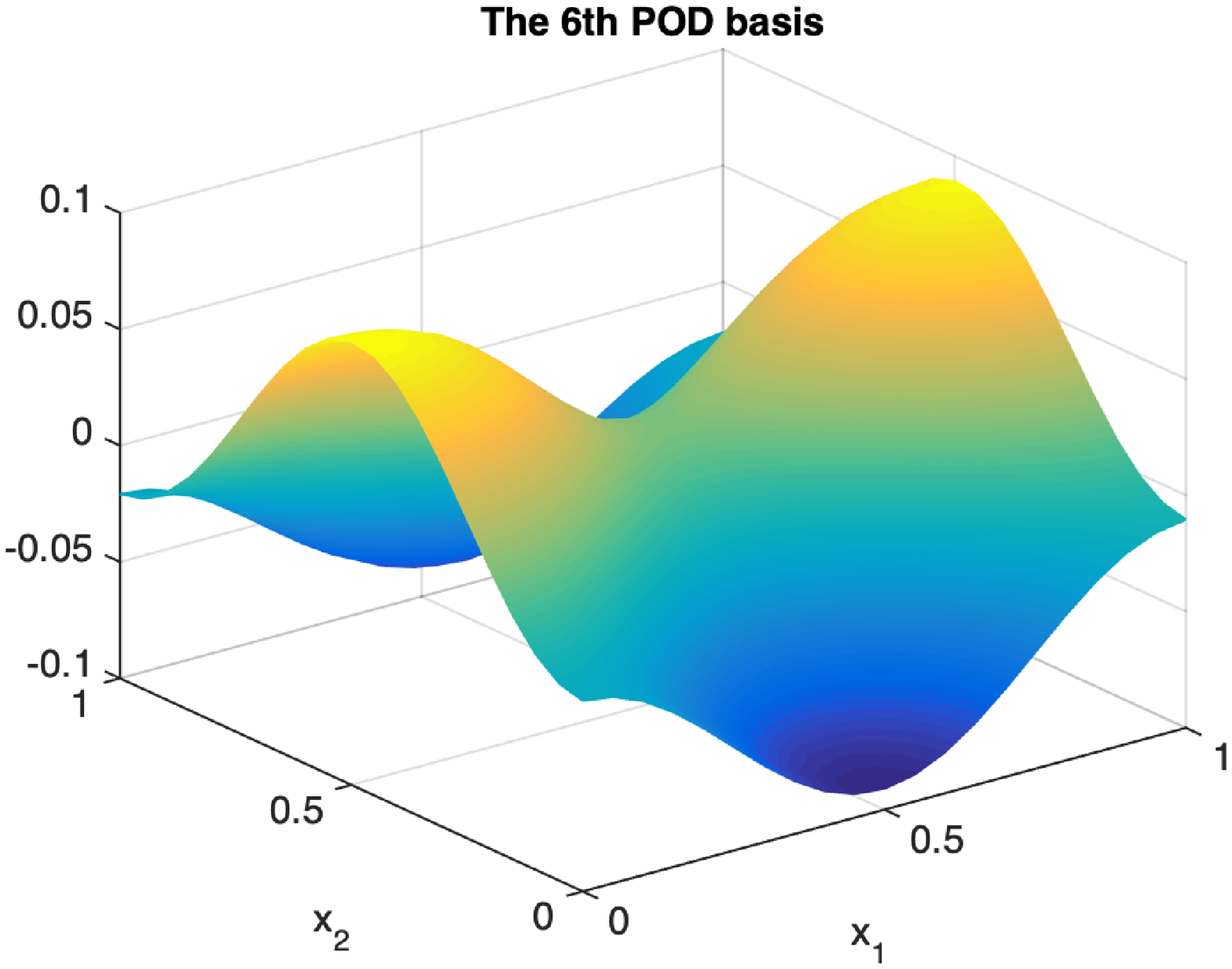}
  \end{overpic}
\end{center}
\caption{ The first six POD basis functions.}\label{RBbases}
\end{figure}

To generate the snapshot set, we sample 100 uniformly distributed parameter location of $\Theta$. Among $100$ time steps, $n_t=50$ ones are included in the database, and this results in an ensemble of $5000$ snapshots. The domain $\Omega$ is discretized to $21\times 21$ regular grid.  The decay of the singular values of the snapshot matrix indicates the approximation power of the POD basis to the information contained in the snapshots. The first 30 singular values are shown in Figure \ref{sigvale}. The singular values decay exponentially, and this indicates that the information contained in the snapshots can be well captured by a few POD basis functions. Through the POD, a set of $p=6$ orthogonal bases are extracted as the RB functions with a tolerance $\epsilon_{pod}=99.99\%$, shown in Figure \ref{RBbases}.

To validate the accuracy of the POD-DSRBF, the RB solutions are recovered at 3 time points for $M=400$ different values of $\theta$. We define the relative  error $\epsilon_a$ and the projection error $\epsilon_p$ as 
\begin{eqnarray*}
\epsilon_a =\frac{1}{M}\sum^M_{i=1} \frac{\|u_h(t;\theta_i)-\tilde{u}_p(t;\theta_i)\|}{\|u_h(t;\theta_i)\|},
\end{eqnarray*}
\begin{eqnarray*}
\epsilon_p =\frac{1}{M}\sum^M_{i=1} \frac{\|u_h(t;\theta_i)-\mb{U}_p\mb{U}_p^Tu_h(t;\theta_i)\|}{\|u_h(t;\theta_i)\|},
\end{eqnarray*}
respectively. We also compute the coefficient learning error $\epsilon_c$
\begin{eqnarray*}
\epsilon_c =\frac{1}{M}\sum^M_{i=1} \frac{\|a(t;\theta_i)-\tilde{a}_p(t;\theta_i)\|}{\|u_h(t;\theta_i)\|}.
\end{eqnarray*}

\begin{figure}
\begin{center}
  \begin{overpic}[width=0.5\textwidth,trim=20 0 20 15, clip=true,tics=10]{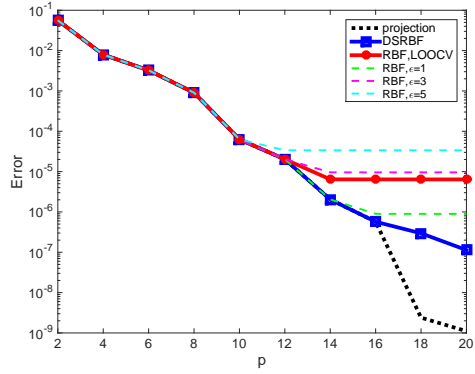}
    \end{overpic}
\end{center}
\caption{The approximation error $\epsilon_a$ as functions of the number of POD basis,  $p$, obtained using the training sets of size $N=200$ and various methods, namely the DSRBF and RBF.  }\label{errorP}
\end{figure}

\begin{figure}
\begin{center}
  \begin{overpic}[width=0.45\textwidth,trim=20 0 20 15, clip=true,tics=10]{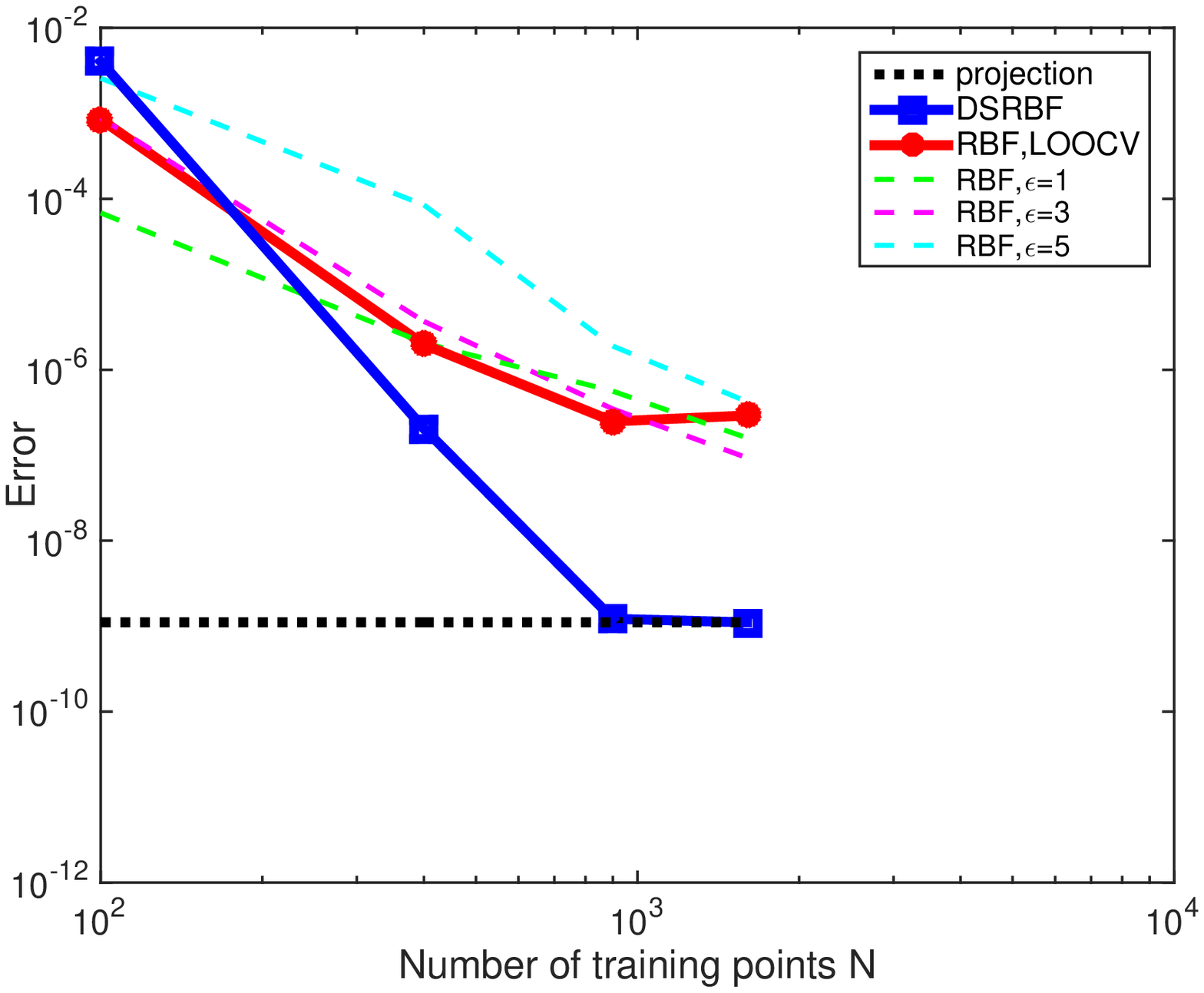}
  \end{overpic}
    \begin{overpic}[width=0.45\textwidth,trim= 20 0 20 15, clip=true,tics=10]{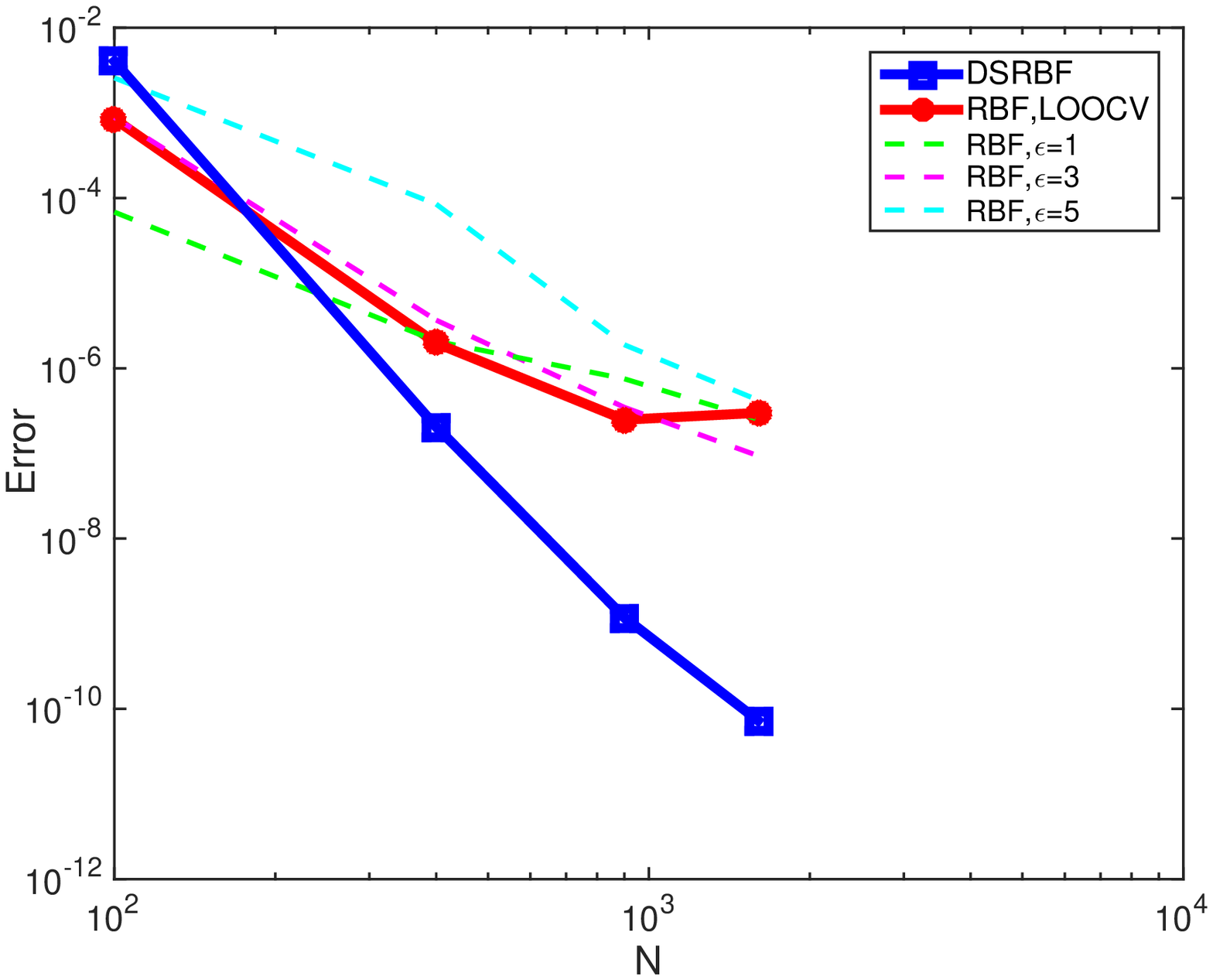}
      \end{overpic}
\end{center}
\caption{The approximation error $\epsilon_a$ (left) and $\epsilon_c$ (right) as functions of the size of training sets,  $N$, obtained using DSRBF and RBF, and $p=20$.}\label{errorQ}
\end{figure}

\begin{figure}
\begin{center}
    \begin{overpic}[width=0.45\textwidth,trim= 20 0 20 15, clip=true,tics=10]{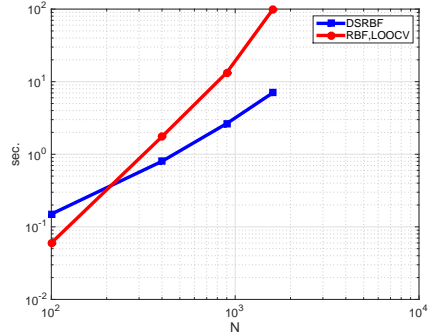}
  \end{overpic}
\end{center}
\caption{CPU times required for the stochastic LOOCV for DSRBF and the LOOCV for RBF.}\label{cpuQ}
\end{figure}

Figure \ref{errorP} illustrates the evolution of the approximation error $\epsilon_a$ and the projection error $\epsilon_p$ with respect to the number of POD basis, $p$, obtained by DSRBF and the classical RBF.  We can observe that both DSRBF and RBF enjoy fast decay with respect to the number of $p$ and resemble the projection  error.  As the number of $p$ increases, classical RBF saturates quickly and while DSRBF can continue to decrease and saturates at a lower level.  The saturation is because the coefficient learning error is dominant over the projection error when the number of $p$ is large.  When more training data is available, the saturation level can be further reduced. Figure \ref{errorQ} shows the approximation error $\epsilon_a$ and the coefficient learning error $\epsilon_c$ with respect to the size of  training data  for a fixed number of $p=20$.  As we increase the number of  training set, both DSRBF and RBF get more accurate results.   Moreover, the DSRBF outperform the classical constant LOOCV-optimal shape parameter formulations. Figure \ref{cpuQ} shows the CPU times for finding LOOCV-optimal shape parameters and the stochastic-LOOCV algorithm respectively.  It can be seen that DSRBF took less than one-tenth of the time of the classical LOOCV algorithm. 
This indicates that the DSRBF can effectively improve the accuracy of the results.

\subsubsection{Solution to the inverse problem} 

In this section, it will be shown how the RB model can be used for accelerating the EKI for solving the model inverse problem, i.e., Example 1.  To better present the results, we shall perform the following two types of approaches: the conventional EKI or direct EKI based on the full model evaluation and RB-EKI based on the POD-DSRBF surrogate model. 

In order to assess the performance of both methods, we define the relative error and data misfit as
$$e_{\theta} = \frac{\|\bar{\theta}-\theta^{\dag}\|}{\|\theta^{\dag}\|},$$ and $$E_{\theta}=\|\Gamma^{-1/2}(y_{obs}-f(\bar{\theta}))\|.$$ Here, $\bar{\theta}$ are $\theta^{\dag}$ are the ensemble mean approximates and exact solutions respectively.  For simplicity, we consider a diagonal measurement error covariance $\Gamma = \delta^2 I$. We first set the exact parameter  $\theta=(0.2,0.7)$. For our EKI methods, we choose $N_e=100$ ensemble members.  

\begin{figure}
\begin{center}
    \begin{overpic}[width=0.45\textwidth,trim= 20 0 20 15, clip=true,tics=10]{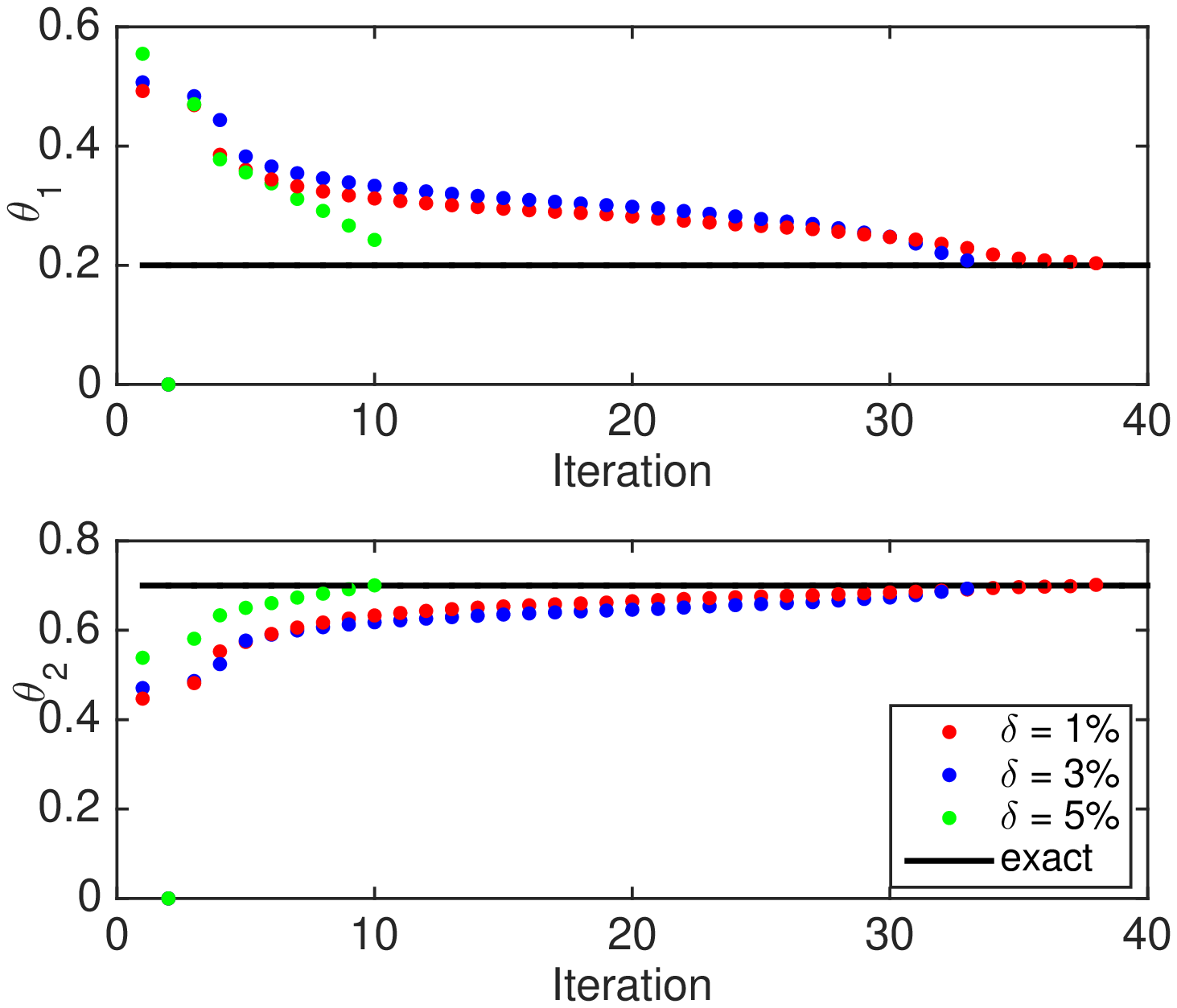}
  \end{overpic}
    \begin{overpic}[width=0.45\textwidth,trim= 20 0 20 15, clip=true,tics=10]{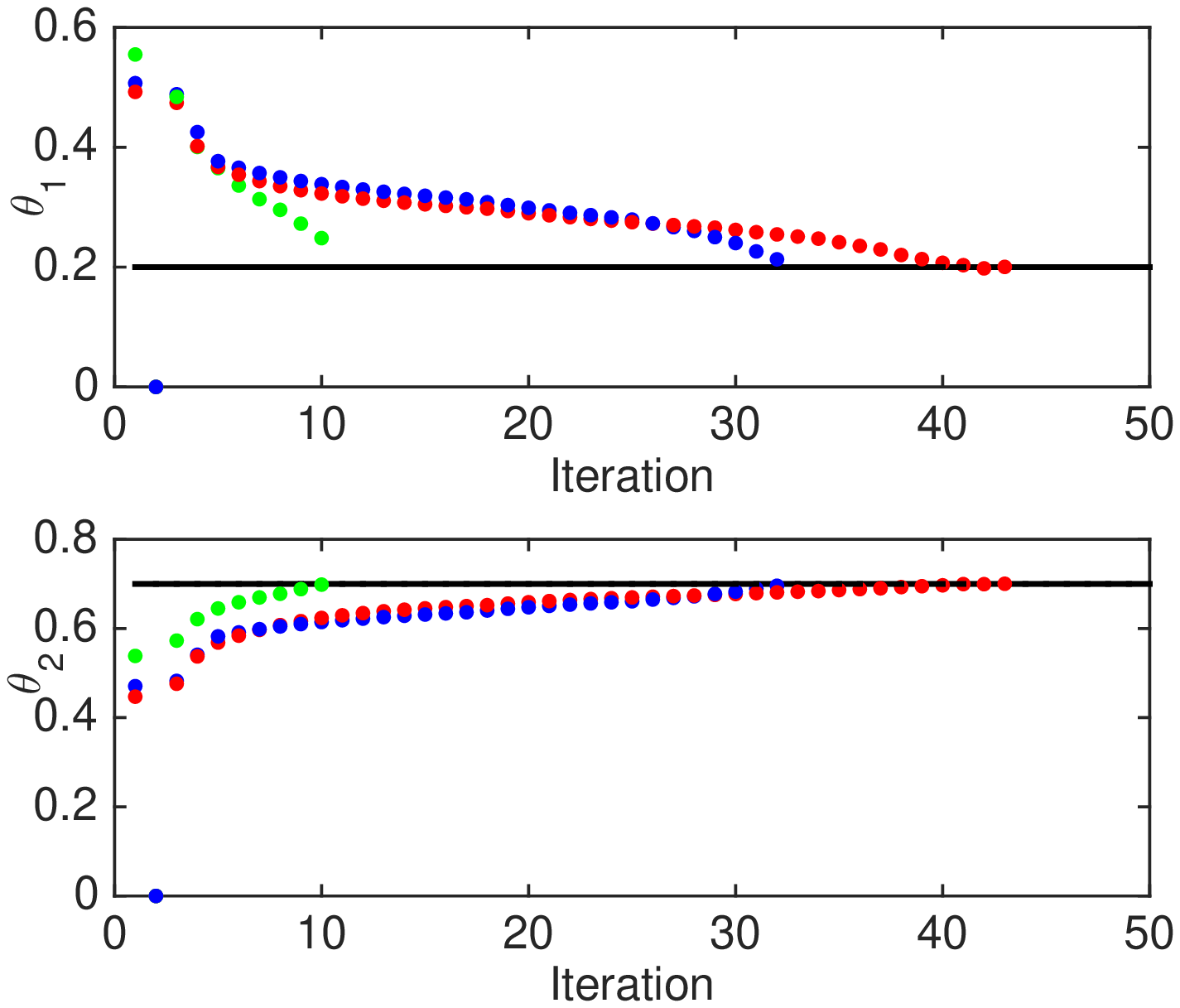}
  \end{overpic}
\end{center}
\caption{The numerical results obtained using full model (left) and RB model (right) with various noise added into the data.}\label{usolut_1}
\end{figure}

\begin{figure}
\begin{center}
  \begin{overpic}[width=0.45\textwidth,trim=20 0 20 15, clip=true,tics=10]{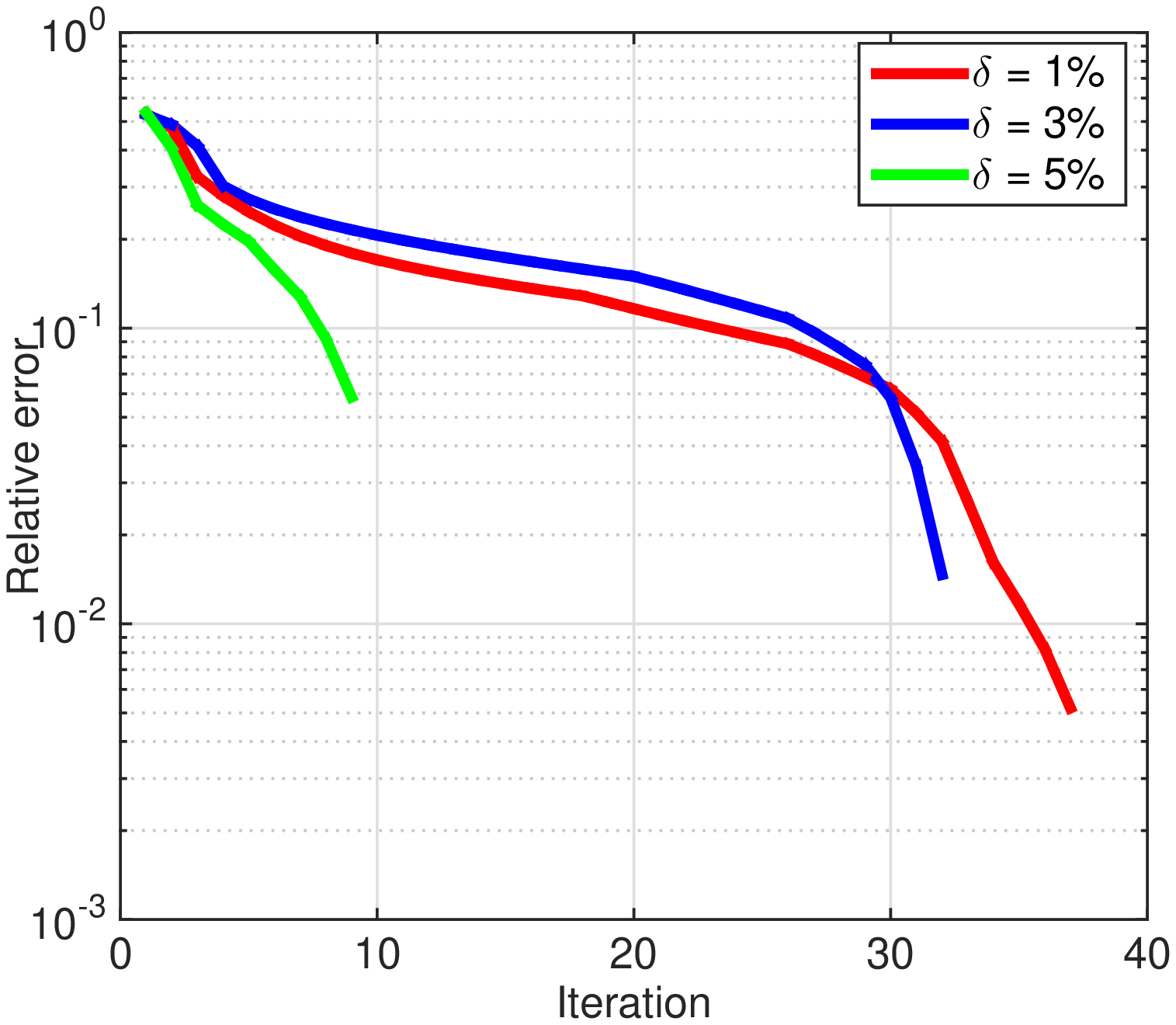}
  \end{overpic}
   \begin{overpic}[width=0.45\textwidth,trim= 20 0 20 15, clip=true,tics=10]{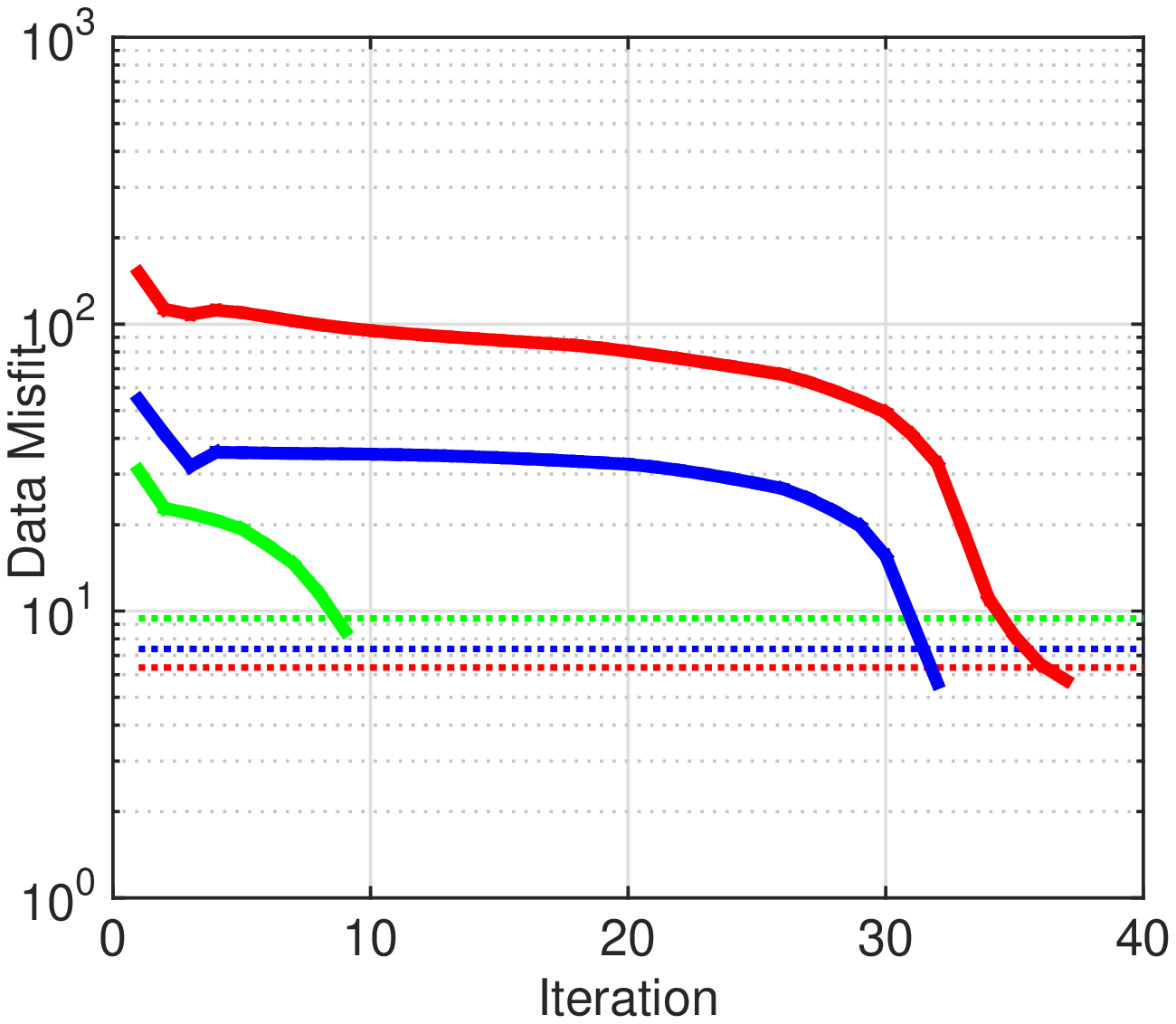}
  \end{overpic}
     \begin{overpic}[width=0.45\textwidth,trim=20 0 20 15, clip=true,tics=10]{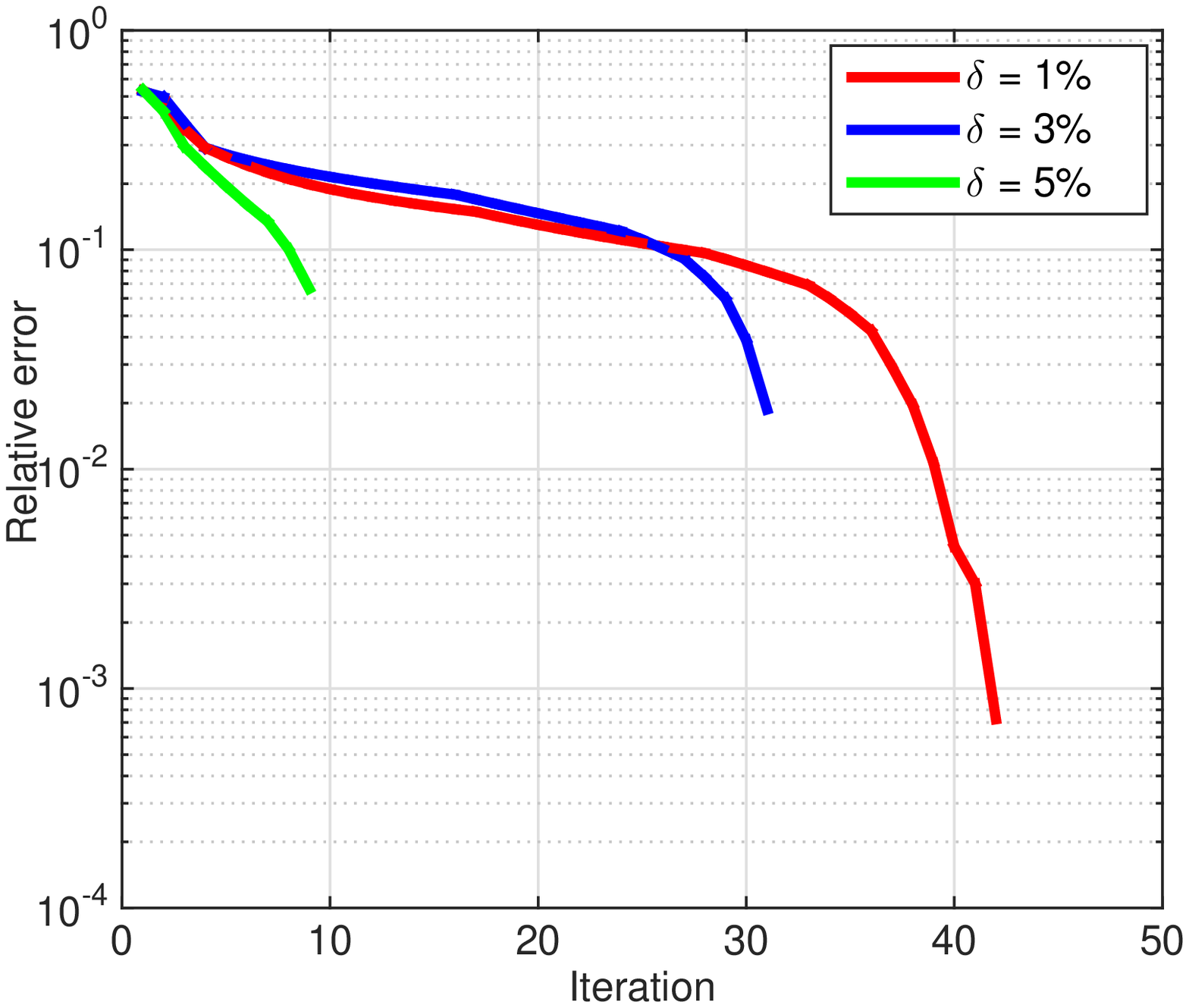}
    \end{overpic}
    \begin{overpic}[width=0.45\textwidth,trim= 20 0 20 15, clip=true,tics=10]{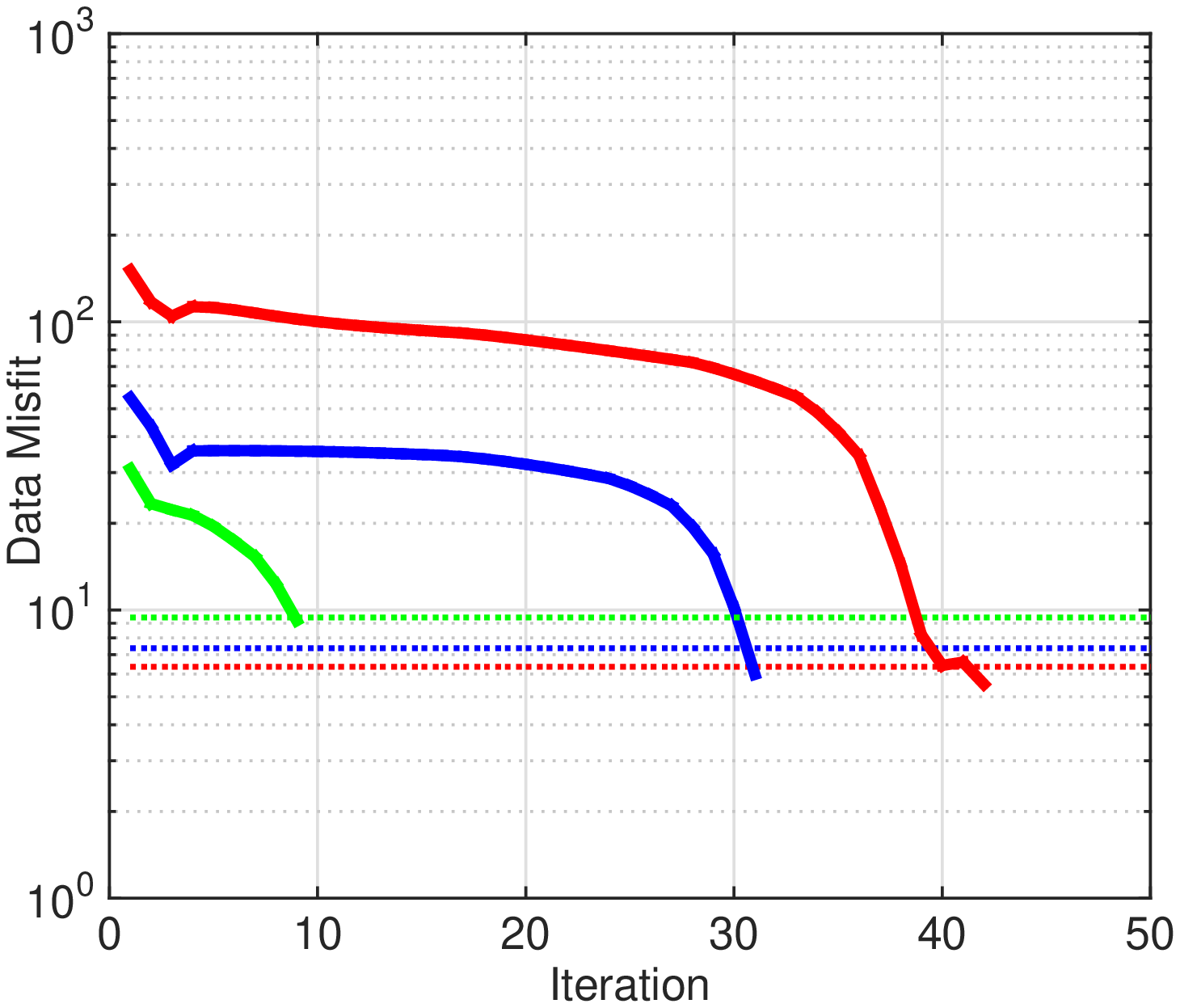}
  \end{overpic}
\end{center}
\caption{The relative errors (left) and data misfits (right) arising from full model(top) and RB model (bottom), respectively. The horizontal line marking the stopping point defined by \eref{dprule}.}\label{unume}
\end{figure}

The exact and numerical results for the source location $\theta$, obtained using full model and RB model with various levels of noise added into the data, namely, $\delta={1\%, 3\%, 5\%}$, are presented in Figure \ref{usolut_1}.  The corresponding  relative error $e_{\theta}$ at each iteration is displayed on the left side of Figure \ref{unume}. From these figures, we can see that though the distance between the initial ensemble mean and the truth is large, the approximation converges towards the exact solution as the iteration progresses.  In particular, the final ensemble mean approximates the truth very well, especially for small noise deviation.  The data misfit along with a line marking the stopping point defined by \eref{dprule} is shown on the right side of Figure \ref{unume}. That is, once the data misfit reaches the dashed line, the convergence criterion defined by \eref{dprule} is considered satisfied. From Figures \ref{usolut_1}  and \ref{unume}, we can see that with the stopping criterion \eref{dprule}, the EKI can yield an accurate and stable approximation. In addition, Figure \ref{unume} clearly shows that both the relative error $e_{\theta}$ and the data misfit $E_{\theta}$  decrease as the iteration progresses.  Moreover, the smaller the deviation in the data is, the more iterations are needed, which means that the cost of obtaining a good approximation from the initial ensemble to the truth increases. This is possibly due to the support of the induced posterior distribution shrinks when the standard deviation in the data decreases \cite{Yan+Zhang2017IP}. 
 
\renewcommand{\arraystretch}{1.5}
\begin{table}
  \centering
  \fontsize{8}{10}\selectfont
  \begin{threeparttable}
    \begin{tabular}{c|cccc|cccc}
    \toprule
    \multicolumn{3}{r}{Direct EKI}&\multicolumn{4}{r}{RB-EKI}\cr
    \cmidrule(lr){2-5} \cmidrule(lr){6-9}
    $\delta (\%)$  
    &$\bar{\theta}^1$ & $\bar{\theta}^2$ &$\# $ iter &CPU(s) &$\bar{\theta}^1$ & $\bar{\theta}^2$ &$\#$iter&CPU(s) \cr
    \midrule
    1   & 0.2034   & 0.7017   & 37    &305.11     & 0.2003  & 0.7004  &42  &4.03\cr
    3   & 0.2082   & 0.6931   & 32    &267.74    & 0.2130  & 0.6959  &31   &2.98\cr
    5  & 0.2426    & 0.7005   & 9     & 76.13      & 0.2484  & 0.6985  &9     &0.87\cr
       \bottomrule
      \end{tabular}
    \end{threeparttable}
      \caption{The numerical results given by two different methods.}\label{uso_eg1_s}
\end{table}

Table \ref{uso_eg1_s} summarize the numerical results with various levels of noise $\delta$ in the data obtained using the direct full model and POD-DSRBF model, respectively.  In the tables, $\bar{\theta}$  denote the final iteration reconstruction arising from EKI.  The ensemble mean of  $\theta$ is in excellent agreement with the exact solution with up to $5\%$ noise in the data. However, the estimate becomes less accurate with an increased noise level, see also Figure \ref{unume}.  The numerical results obtained by the two approaches are nearly the same, which is consistent with the discussion from the last section that the POD-DSRBF model is accurate enough to approximate the full model.  The foregoing numerical verifications indicate that the EKI with the full model simulation and POD-DSRBF model could yield practically identical numerical results for the  inverse heat source problem.   However, the computing time required by RB-EKI is only a small fraction of that by direct full model simulation. To see that, consider the case with $\delta = 3\%$, the RB-EKI takes $2.98s$ to complete the simulation, whereas direct EKI takes $267.74s$.  Therefore, the speedup of the RB-EKI approach over the direct EKI approach is dramatic.


Besides the great computational savings, another advantage of using RB-EKI is that it is reusable when there are new measurements coming after constructing the POD-DSRBF model. Table \ref{meanP2} shows the results for four other source locations with $3\%$ noise in the data. The POD-DSRBF model is the same as before. It can be seen that the method can always infer the exact value without performing any additional direct forward simulation as in the direct EKI.  This demonstrates the flexibility of the RB-EKI approach for solving the time fractional heat source identification problems.

  \begin{table} 
   \begin{center}
  \arrayrulewidth=1pt
  \renewcommand{\arraystretch}{1.5}
\begin{tabular}{lllll}
  \hline
  $ \theta$ & $\bar{\theta}^1$ & $\bar{\theta}^2$  & $\#$ iter& CPU(s)  \\
  \hline
(0.8, 0.4)  & 0.7969 & 0.4101 & 30 & 2.89 \\
 (0.8, 0.2) & 0.7572 & 0.1736 & 30 & 2.88 \\
 (0.9, 0.2) & 0.7691 & 0.1032 & 47 & 4.53 \\
 (0.9, 0.1) & 0.7765 & 0.0483 & 46 & 4.43 \\
  \hline
\end{tabular}
 \end{center}
  \caption{The numerical results  with four other source locations given by RB-EKI.}\label{meanP2}
 \end{table} 
 
 \subsubsection{Possibility of recovering the fractional order $\alpha$}
 
In this section, we consider  to recover the fractional order $\alpha$ in the model (\ref{2dsource})  with the source location simultaneously. In this case, the unknown parameter is $\theta=(\theta^1,\theta^2,\alpha)\in \R^3$. For the results presented below, a uniform prior was adopted  for $\theta$ on the parameter space $[0,1]^3$.  To generate data for inversion, the `true' parameter is choose to be $(\theta^1=0.25, \theta^1=0.75, \alpha=0.8)$. In the following tests, the measurements were taken at three different times $t=\{0.25, 0.75,1\}$ using 25 locations on  a uniform $5\times5$ grid covering the domain $\Omega$, thus leading to $75$ measurements. To construct the POD-DSRBF model, we sample the parameter space $[0,1]^3$ at 200 uniformly distributed random points, and add the full-order solutions at 50 time points $t\in \mathcal{T}_{tr}=\{0.02,0.04,\cdots,1\}$ into the database $\mathcal{D}$. These samples will be used as snapshots and training data.  In this example, we use $p=10$ POD basis functions to construct the RB space,  and the truncation number of SVDs for all expansion coefficients are set as $q_k=10$. 

 \begin{figure}
\begin{center}
  \begin{overpic}[width=0.45\textwidth,trim=20 0 20 15, clip=true,tics=10]{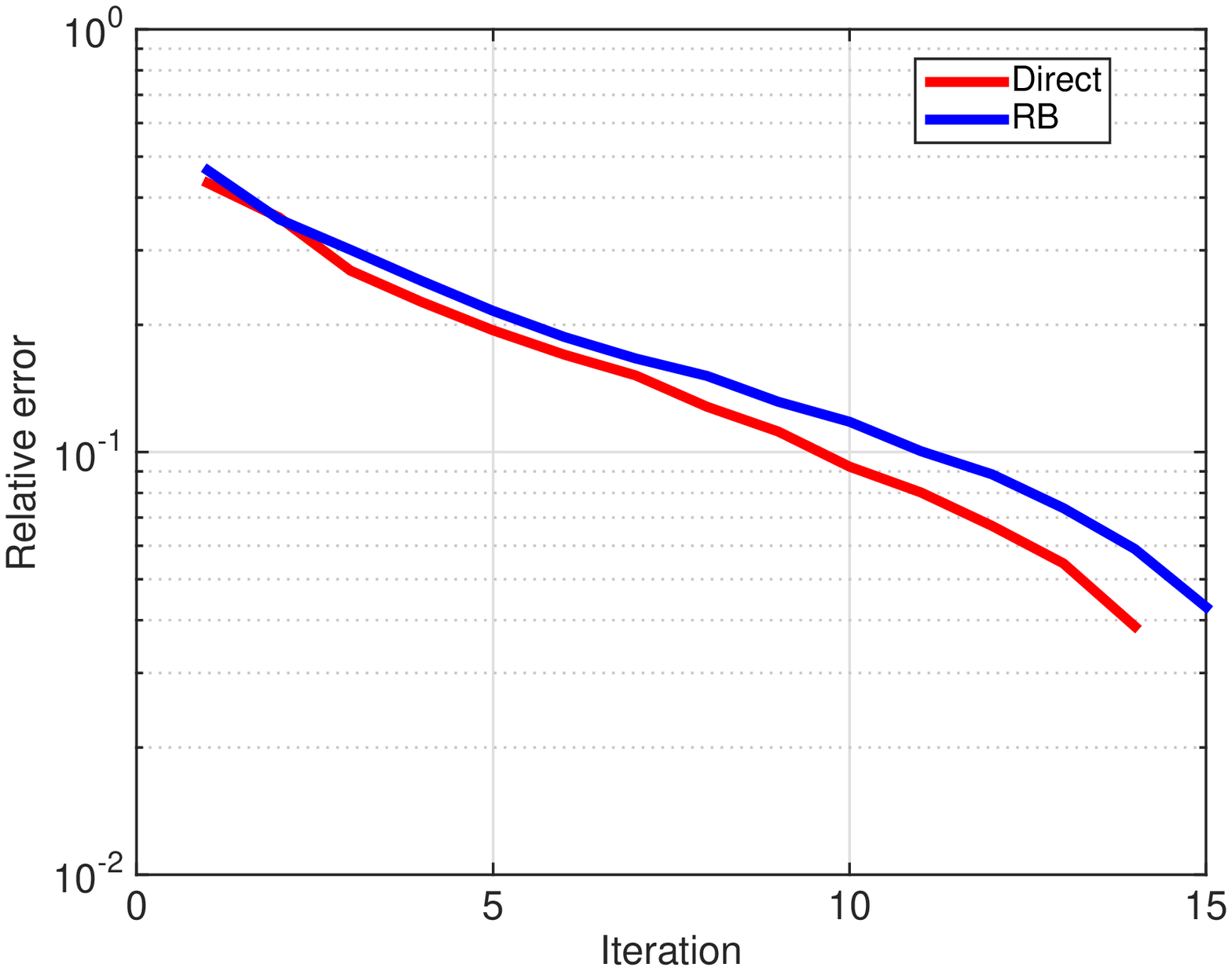}
  \end{overpic}
    \begin{overpic}[width=0.45\textwidth,trim= 20 0 20 15, clip=true,tics=10]{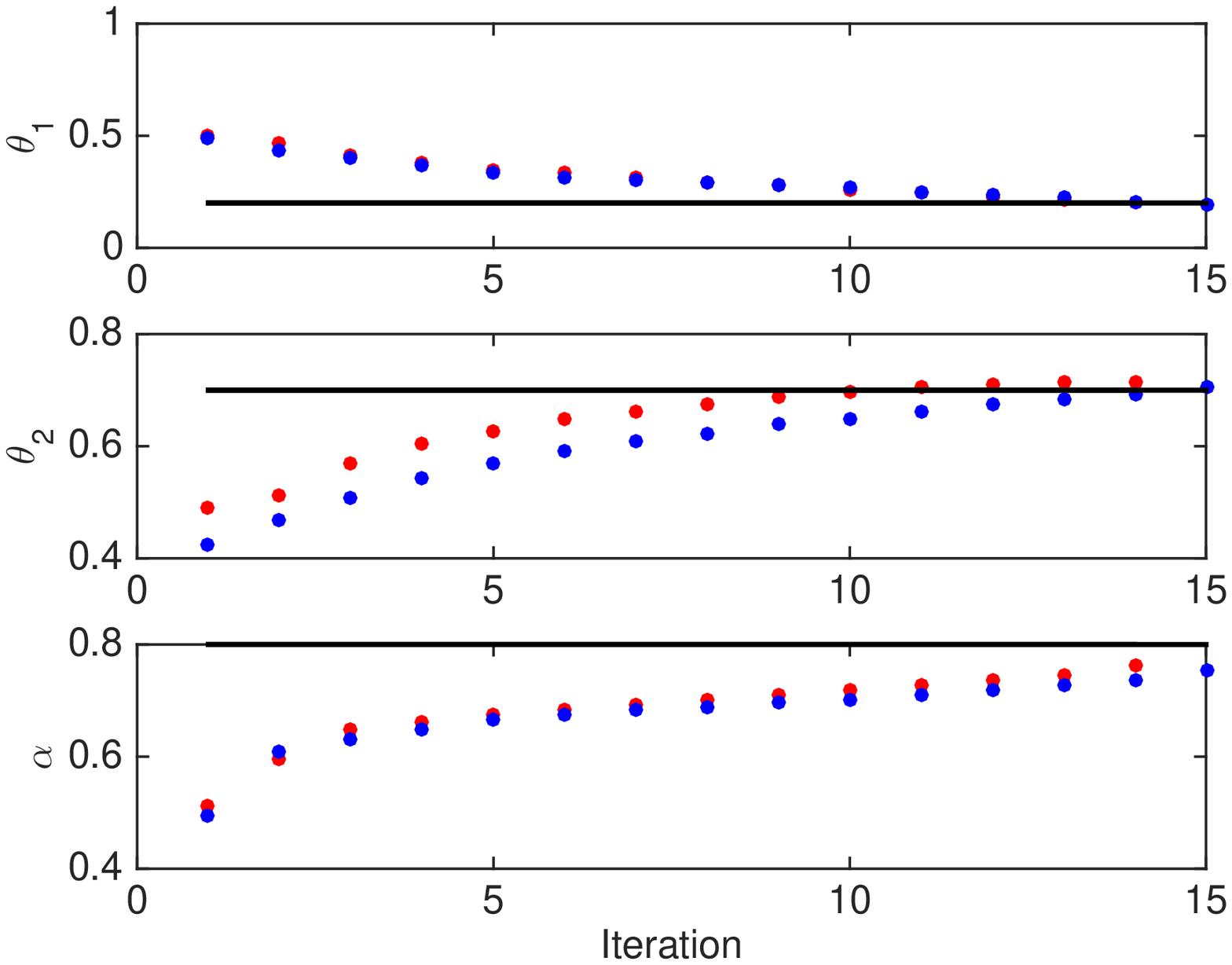}
  \end{overpic}
\end{center}
\caption{The numerical results arising full model and RB model with $5\%$ noise in the data.}\label{u_eg2}
\end{figure}

The numerical results by RB-EKI and the direct EKI with $5\%$ noise in the data are shown in Figure \ref{u_eg2}.  The approximation converges towards the exact solution as the iteration progresses. In particular, the final ensemble mean  is in excellent agreement with the exact one for both RB-EKI and direct EKI. The result by RB-EKI is practically identical to that by the direct EKI, but its computational expense is only a small fraction of that of direct EKI.  The CPU time of evaluating the conventional EKI is about $256.86s$, while the CPU time of RB-EKI is about $2.82s$.  

\subsection{Example 2: Estimating the diffusion coefficient}
In this example, we illustrate the RB-EKI approach on the nonlinear inverse problem of estimating the  diffusion coefficient. Consider the following two dimensional time-fractional PDEs
\begin{eqnarray}\label{2dtfpde}
\begin{array}{rl}
^cD_t^{\alpha}u-\nabla\cdot(\kappa(x) \nabla u(x,t))&=e^{-t}\exp\Big(-\frac{\|x-(0.25,0.75)\|^2}{2\times0.1^2}\Big),\quad \Omega\times [0, 1],\\
\nabla u \cdot \textbf{n}&=0, \quad \mbox {on} \,\partial{\Omega},\\
u(x,0)&=0, \quad  \mbox{in}\, \Omega.
 \end{array}
\end{eqnarray}
The goal is to determine the diffusion coefficient  $\kappa(x)$ from noisy measurements of the $u$-field at a finite set of locations and times. We assume the log-diffusivity field $\log\kappa(x)$ is endowed with a Gaussian process prior, with mean zero and an isotropic  covariance kernel:
\begin{equation*}
K(x_1,x_2)=\sigma^2 \exp\Big(-\frac{\|x_1-x_2\|^2}{2l^2}\Big),
\end{equation*}
for which we choose variance $\sigma^2=1$ and a length scale $l = 0.2$. This prior allows the field to be easily parameterized with a Karhunen-Loeve expansion:
\begin{equation}
\log\kappa(x;\theta) \approx \sum^{d}_{i=1} \theta^i \sqrt{\lambda_i} \phi_i(x),
\end{equation}
where $\lambda_i$ and $\phi_i(x)$ are  the eigenvalues and eigenfunctions, respectively, of the integral operator on $[0,1]^2$ defined by the kernel $K$, and the parameter $\theta^i$ are endowed with independent standard normal priors. These parameters then become the targets of inference.  In the numerical simulation, we truncate the Karhunen-Loeve expansion at $d=9$ modes that preserve $90\%$ energy of the prior distribution.  The true solution is directly drawn from the prior distribution.  The measurement sensors of $u$ are using 49 locations on uniform $7 \times 7$ grid covering the domain $\Omega$.  Similar to Example 1, at each sensor location, three measurements are taken at  time $t=\{0.25,0.75,1\}$. The observational errors are taken to be additive and Gaussian:
\begin{equation*}
y_j = u(x_j,t_j;\theta) +\xi_j,
\end{equation*}
with $\xi_j \sim N(0,\delta^2)$.

\begin{figure}
\begin{center}
  \begin{overpic}[width=0.45\textwidth,trim=20 0 20 15, clip=true,tics=10]{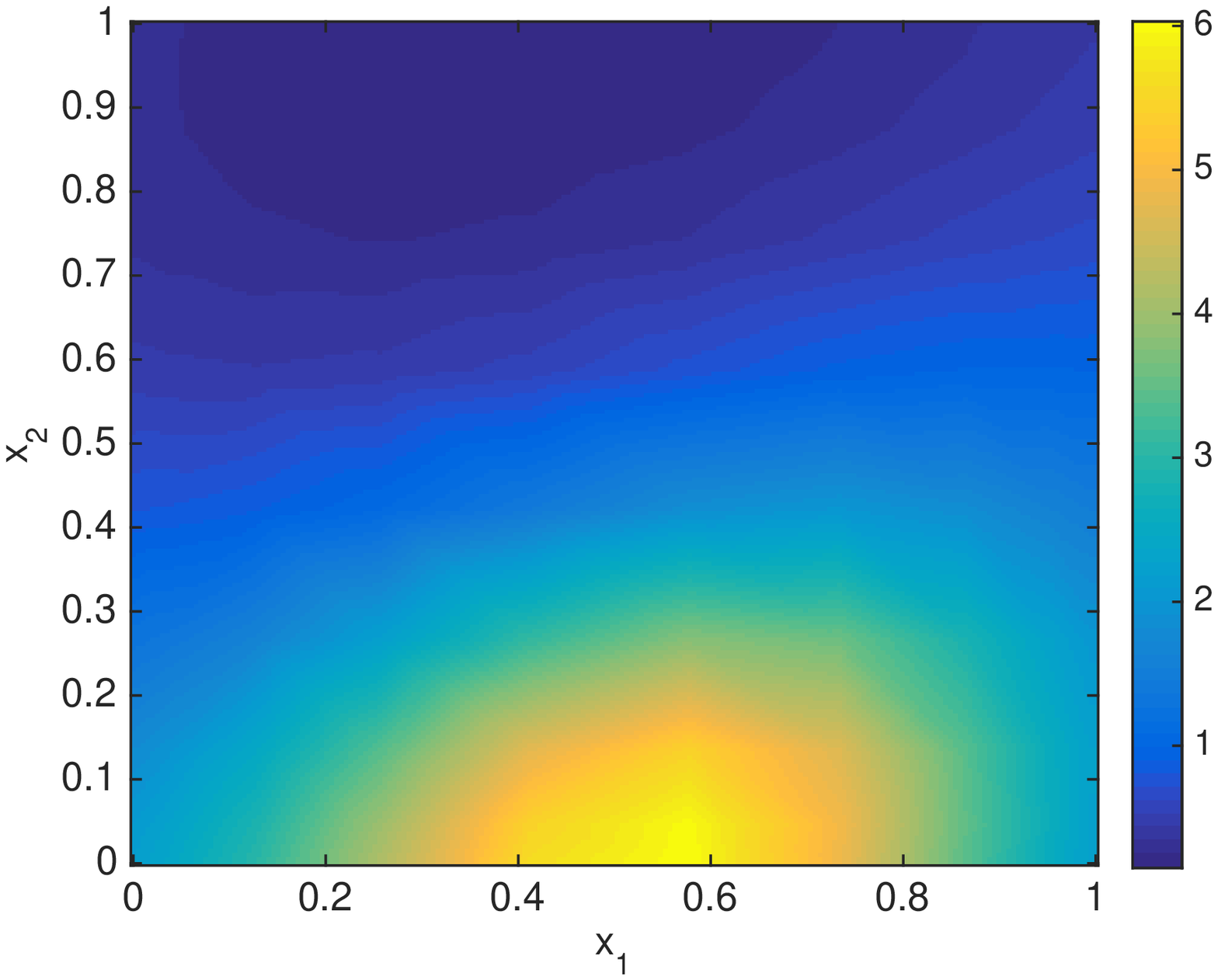}
  \end{overpic}
    \begin{overpic}[width=0.45\textwidth,trim= 20 0 20 15, clip=true,tics=10]{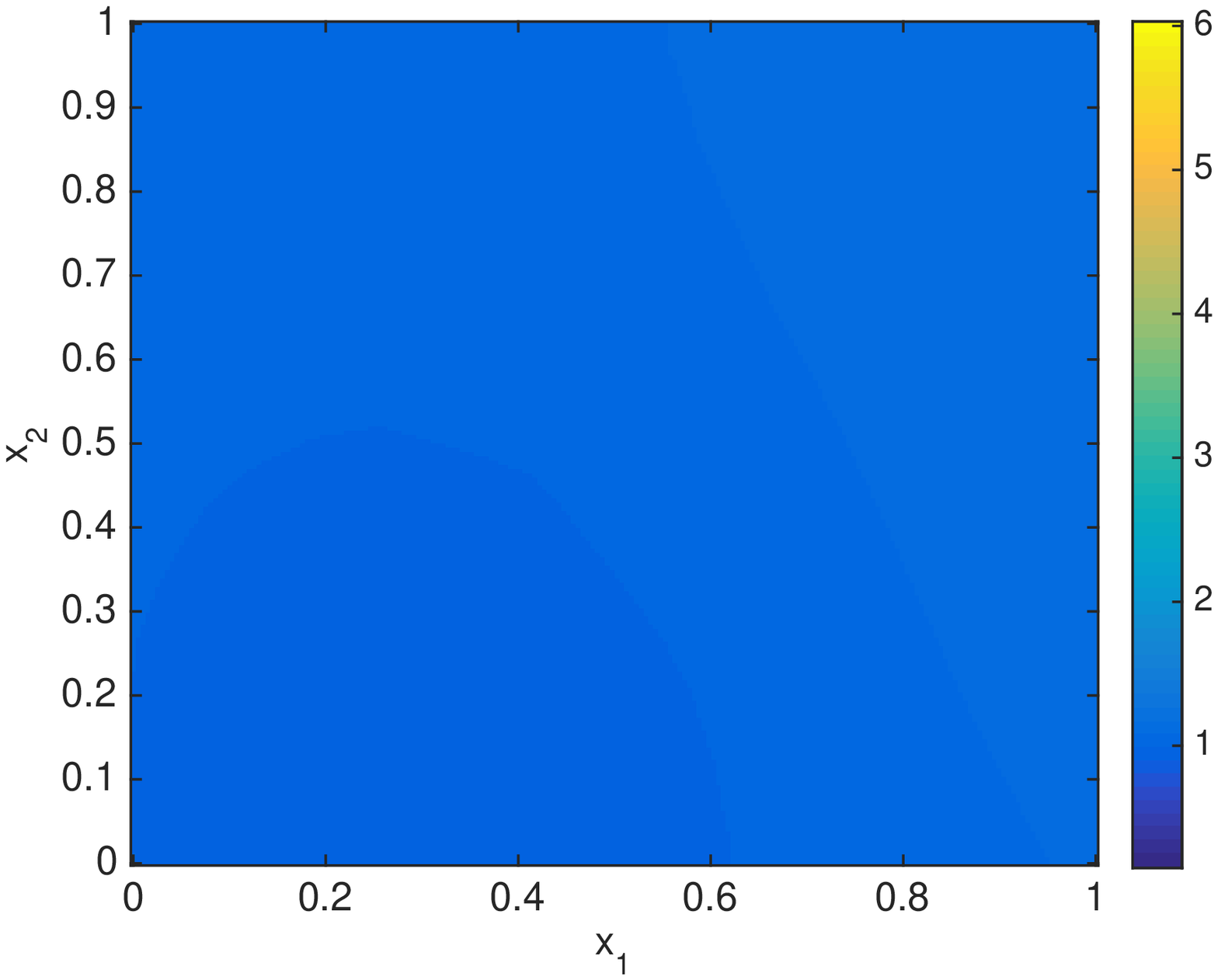}
  \end{overpic}
\end{center}
\caption{ Left: the true diffusion coefficients used for generating the synthetic data set. Right: the initial ensemble mean.}\label{uexact_eg3}
\end{figure}

\begin{figure}
\begin{center}
  \begin{overpic}[width=0.32\textwidth,trim=20 0 20 15, clip=true,tics=10]{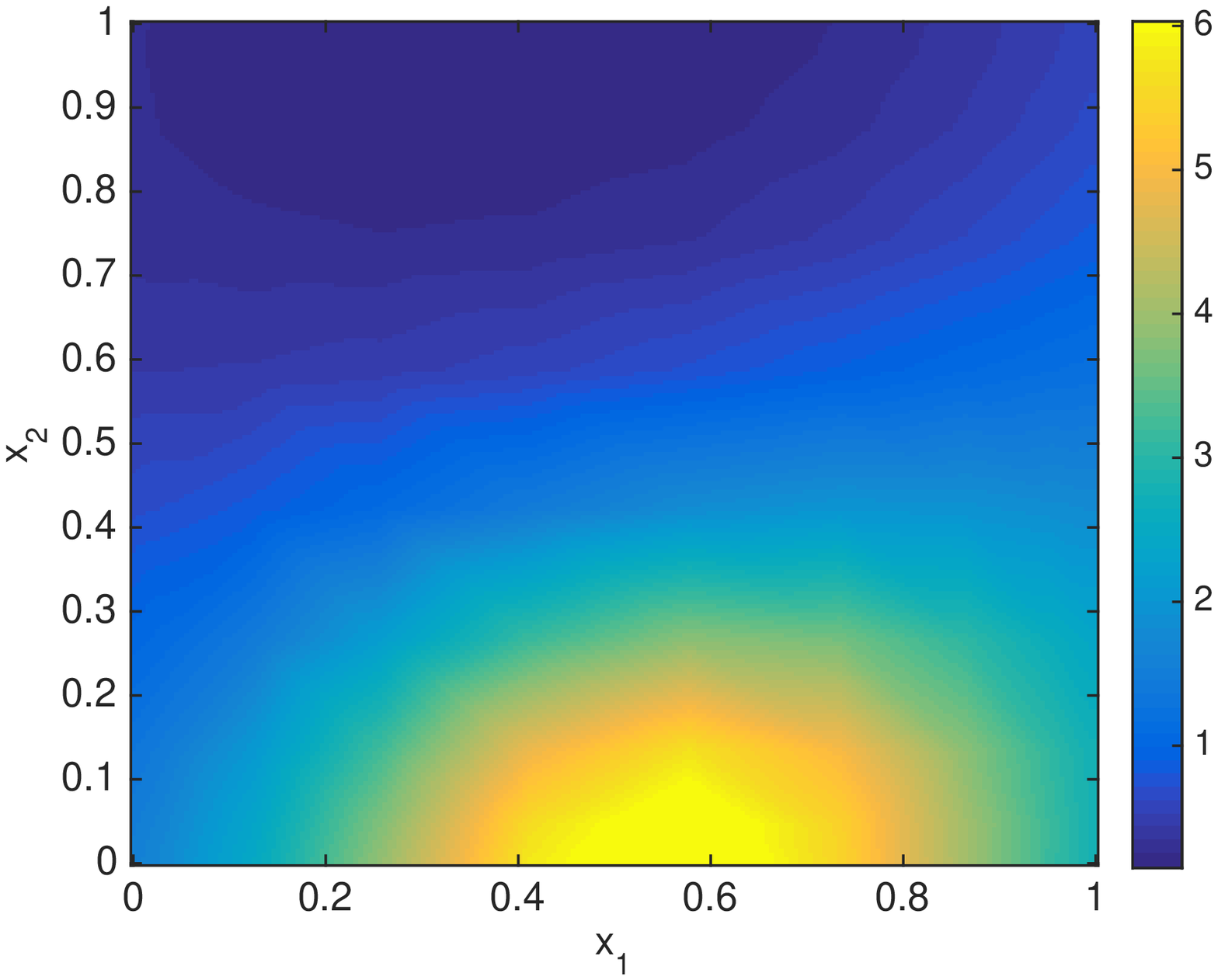}
  \end{overpic}
    \begin{overpic}[width=0.32\textwidth,trim= 20 0 20 15, clip=true,tics=10]{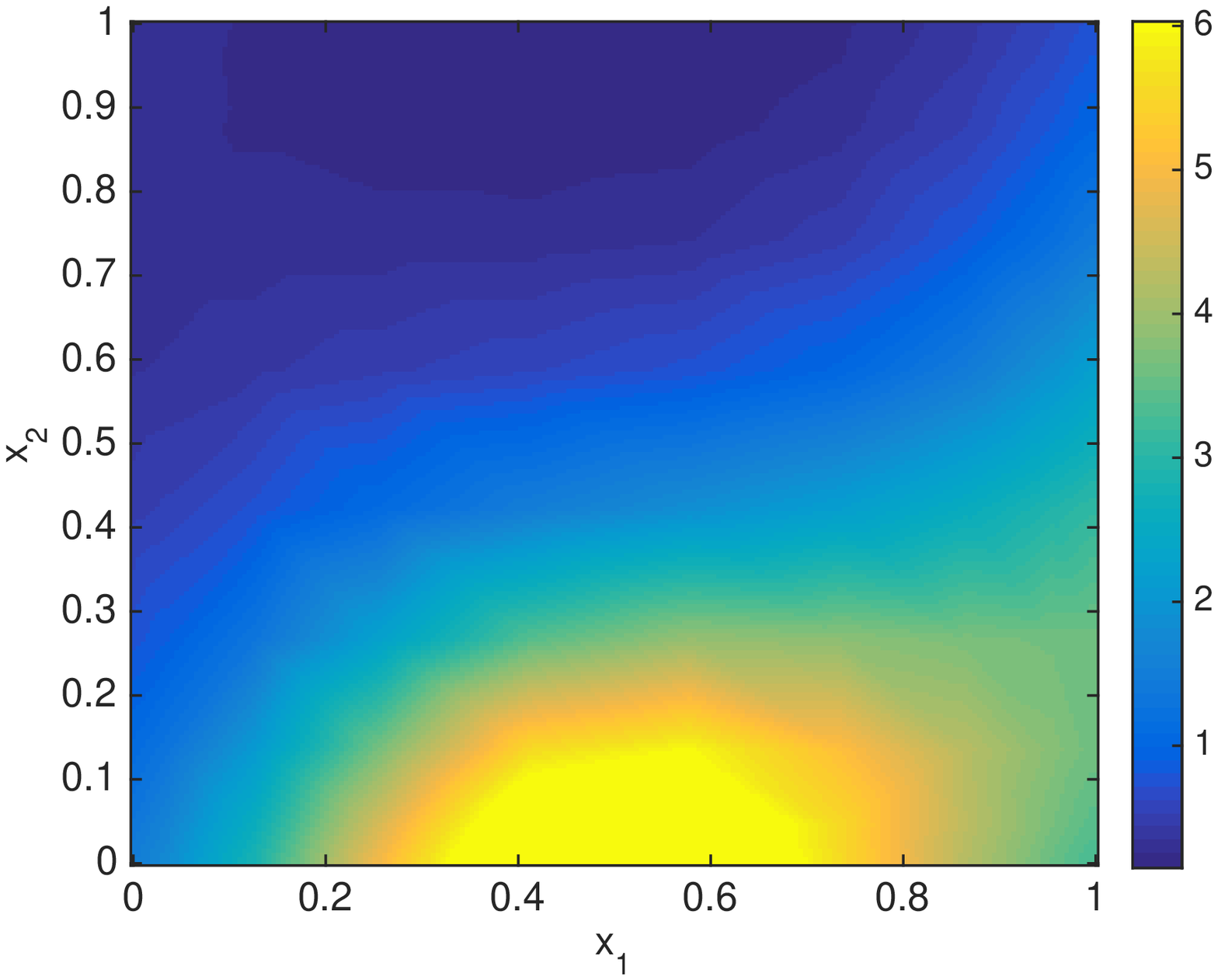}
  \end{overpic}
     \begin{overpic}[width=0.32\textwidth,trim= 20 0 20 15, clip=true,tics=10]{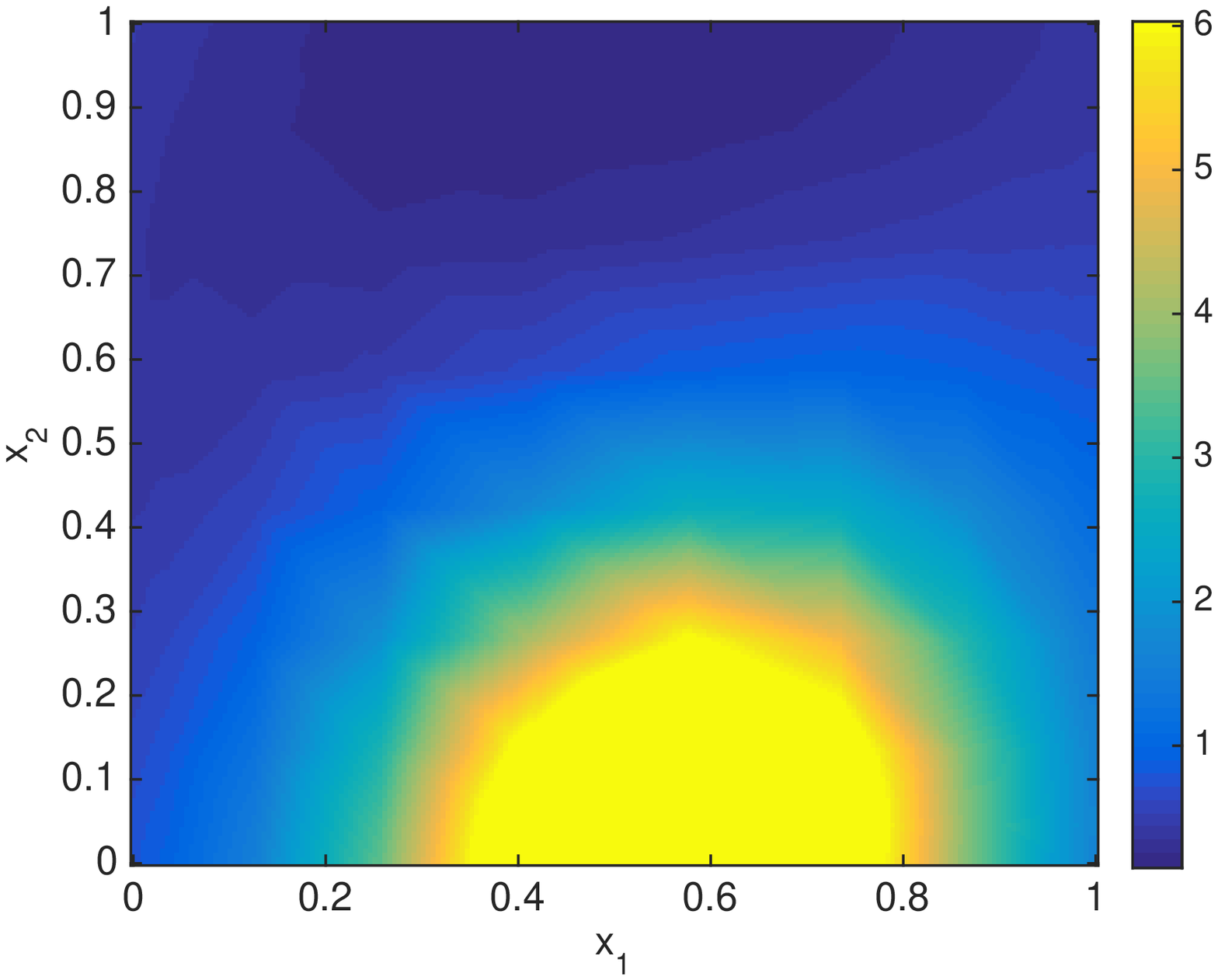}
  \end{overpic}
\end{center}
\caption{Numerical results for the final iteration by RB EKI, using $N_e = 200$ and various levels of noise in the data (from left to right): $1\%, \, 3\%$ and $5\%$.}\label{u_eg3}
\end{figure}

To construct the POD-DSRBF model, we sample the parameter space  at 500 random points, and the full-order database is constructed in the same way as Example 1.  All full-order samples will be used as snapshots and training data.  It is observed that $p=10$ POD basis functions already portray accurately the solution manifold, with the total energy $\epsilon_{pod}=99.99\%$. The true diffusion coefficient and the initial ensemble mean of the EKI are shown in Figure \ref{uexact_eg3}. The final iteration reconstruction by RB-EKI for Example 2 with various levels of noise in the data is shown in Figure \ref{u_eg3}.  With up to $5\%$ noise in the data, the ensemble mean is still in good agreement with the exact one.  Taking into consideration the ill-posedness of the problems, the results presented here are quite satisfactory.  Direct EKI simulation is not conducted in this example since it is expected to much more expensive.

\begin{figure}
\begin{center}
  \begin{overpic}[width=0.45\textwidth,trim=20 0 20 15, clip=true,tics=10]{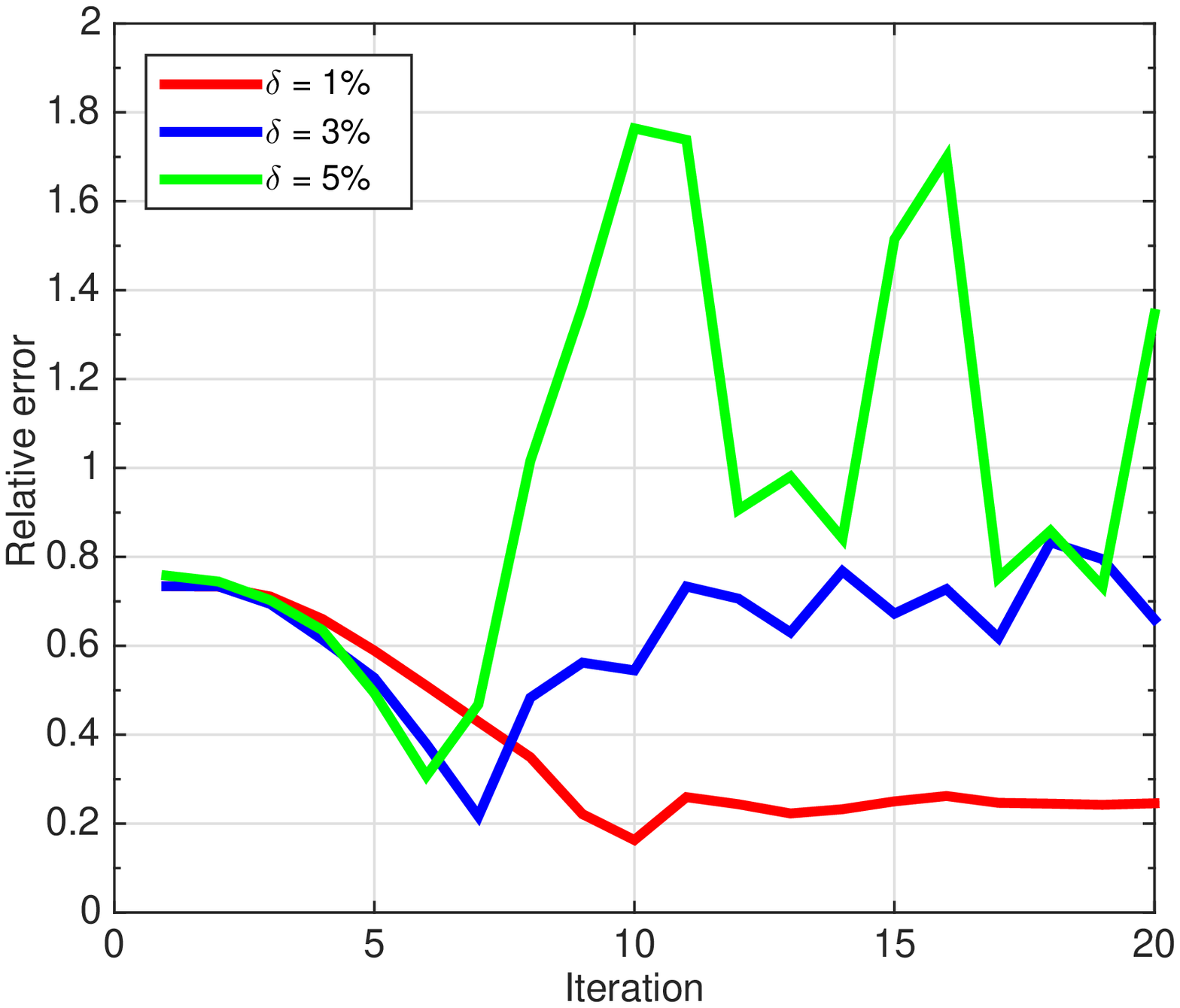}
  \end{overpic}
    \begin{overpic}[width=0.45\textwidth,trim= 20 0 20 15, clip=true,tics=10]{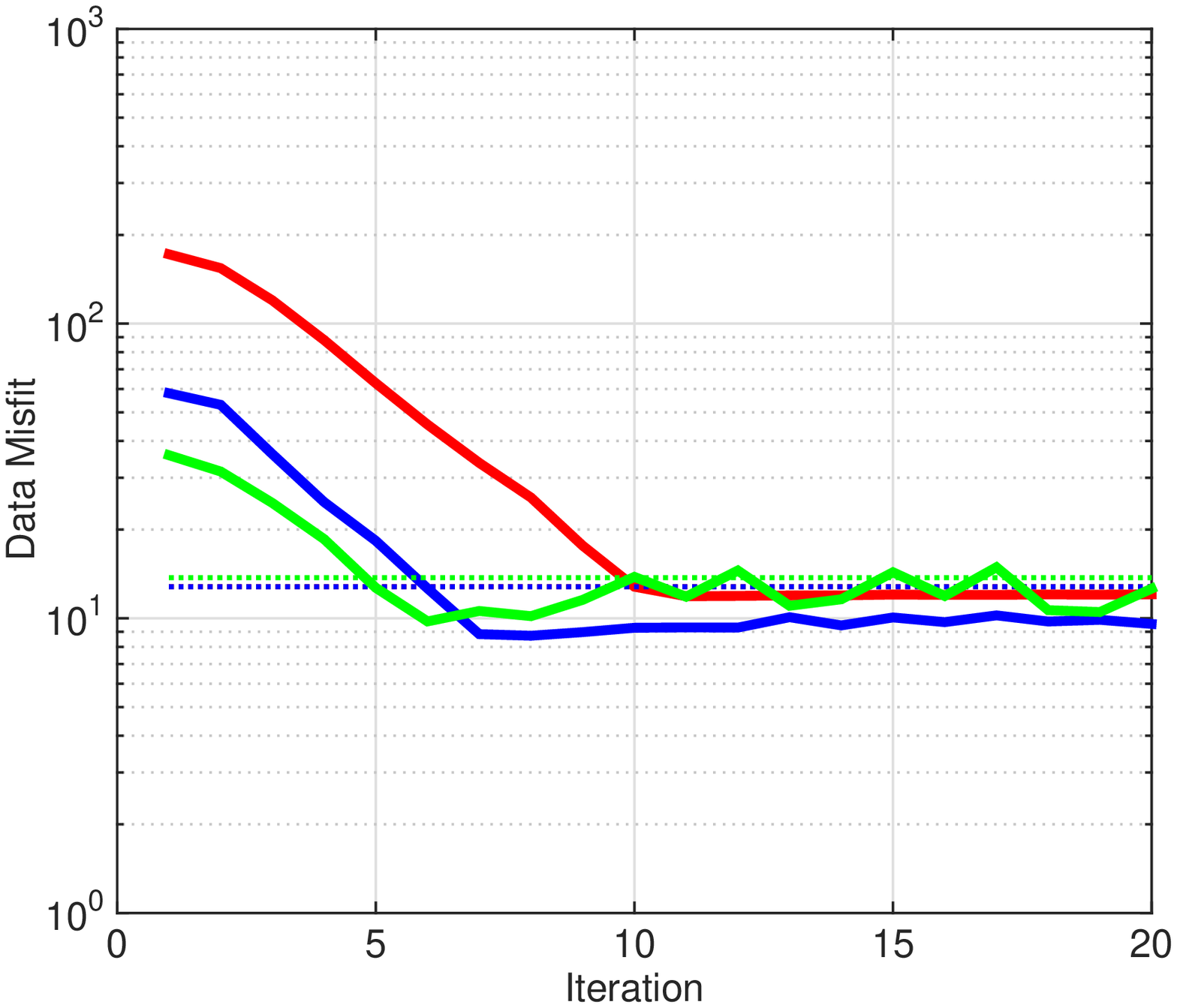}
  \end{overpic}
\end{center}
\caption{ Relative errors (left) and data misfits (right) of RB EKI. The horizontal line marking the stopping point defined by \eref{dprule}.}\label{err_eg3}
\end{figure}

In order to study the convergence of the results, we also plot the relative errors and data misfits of RB-EKI at  each iteration. The numerical results together with a horizontal line marking the stopping points defined by the \eref{dprule} are  displayed in Figure \ref{err_eg3}.  The accuracy error $e_{\theta}$ of the RB-EKI exhibits the typical `semiconvergence' phenomenon: the approximation solution converges towards the exact solution up to a certain iteration number, and beyond this point, it deviates from the exact solution. Therefore, an appropriate stopping rule is crucial to obtain an accurate and stable approximation solution.  The right side of Figure \ref{err_eg3} shows that the data misfits curve is relatively flat over a wide range of iteration numbers, and thus the rule \eref{dprule}  requires an accurate estimation of the amount of data noise, in order that the RB-EKI could yield a good approximation. 
However, this is not always available in practical situations, and heuristic approaches, such as the generalized cross-validation and the L-curve criterion, may be employed.  We leave this to our next exploration.

\section{Summary} \label{sec:summary}

In this paper we developed an efficient non-intrusive RB ensemble Kalman inversion approach to solving the time fractional diffusion inverse problems. In the framework of this method, a reduced order model is constructed by POD-DSRBF.  The identification of the optimal shape parameter of the DSRBF is performed during the offline stage through a stochastic LOOCV algorithm.  Then the prediction steps of the EKI are calculated from a large number of realizations generated by the reduced order model with virtually no additional computational cost.  Thus, the computational cost can be significantly reduced. With the accuracy and the efficiency validated by numerical examples, the proposed RB-EKI algorithm is shown to be a powerful tool for solving time fractional diffusion inverse problems. 

Ongoing works will extend the techniques to more complex fractional order model inverse problems, e.g. time-space fractional diffusion equation.  Note that the reduced models under investigation are global in that they are valid over the prior parametric domain, which underlies the challenge of the high dimensionality of the parameter space as well as the intrinsic properties of the forward model.  Strategies for building multi-fidelity reduced models adaptively would be very interesting.  

\bibliographystyle{unsrt}
\bibliography{mEnPCK}

\begin{thebibliography}{10}

\bibitem{Podlubny1999}
I.~Podlubny.
\newblock {\em Fractional differential equations}, volume 198 of {\em
  Mathematics in Science and Engineering}.
\newblock Academic Press Inc., San Diego, CA, 1999.

\bibitem{SY2011}
K.~{Sakamoto} and M.~{Yamamoto}.
\newblock Initial value/boundary value problems for diffusion-wave equations
  and applications to some inverse problems.
\newblock {\em Journal of Mathematical Analysis and Applications},
  382(1):426--447, 2011.

\bibitem{Lin+Xu2007}
Y.~Lin and C.~Xu.
\newblock Finite difference/spectral approximations for the time-fractional
  diffusion equation.
\newblock {\em Journal of Computational Physics}, 225(2):1533--1552, 2007.

\bibitem{FCY2013JCP}
Z.~J. Fu, W.~Chen, and H.~T. Yang.
\newblock Boundary particle method for laplace transformed time fractional
  diffusion equations.
\newblock {\em Journal of Computational Physics}, 235:52--66, 2013.

\bibitem{JR2015}
B.~T. {Jin} and W.~{Rundell}.
\newblock A tutorial on inverse problems for anomalous diffusion processes.
\newblock {\em Inverse Problems}, 31(3):35003, 2015.

\bibitem{CNYY2009}
J.~{Cheng}, J.~{Nakagawa}, M.~{Yamamoto}, and T.~{Yamazaki}.
\newblock Uniqueness in an inverse problem for a one-dimensional fractional
  diffusion equation.
\newblock {\em Inverse Problems}, 25(11):115002, 2009.

\bibitem{JPY2017}
J.~X. {Jia}, J.~G. {Peng}, and J.~Q. {Yang}.
\newblock Harnack's inequality for a space-time fractional diffusion equation
  and applications to an inverse source problem.
\newblock {\em Journal of Differential Equations}, 262(8):4415--4450, 2017.

\bibitem{LZJY2013}
G.~S. {Li}, D.~L. {Zhang}, X.~Z. {Jia}, and M.~{Yamamoto}.
\newblock Simultaneous inversion for the space-dependent diffusion coefficient
  and the fractional order in the time- fractional diffusion equation.
\newblock {\em Inverse Problems}, 29(6):65014, 2013.

\bibitem{LYY2015}
J.J. Liu, M.~Yamamoto, and L.~Yan.
\newblock On the reconstruction of unknown time-dependent boundary sources for
  time fractional diffusion process by distributing measurement.
\newblock {\em Inverse Problems}, 32(1):015009, 2016.

\bibitem{YY2015}
L.~Yan and F.~L. Yang.
\newblock The method of approximate particular solutions for the
  time-fractional diffusion equation with a non-local boundary condition.
\newblock {\em Computers \& Mathematics with Applications}, 70(3):254--264,
  2015.

\bibitem{ZX2011}
Y.~{Zhang} and X.~{Xu}.
\newblock Inverse source problem for a fractional diffusion equation.
\newblock {\em Inverse Problems}, 27(3):35010, 2011.

\bibitem{ZW2010}
G.~H. {Zheng} and T.~{Wei}.
\newblock Two regularization methods for solving a riesz-feller
  space-fractional backward diffusion problem.
\newblock {\em Inverse Problems}, 26(11):115017, 2010.

\bibitem{TS2017}
N.~G. {Trillos} and D.~{Sanz-Alonso}.
\newblock The bayesian formulation and well-posedness of fractional elliptic
  inverse problems.
\newblock {\em Inverse Problems}, 33(6):65006, 2017.

\bibitem{ZJY2018IP}
Y.~X. Zhang, J.~Jia, and L.~Yan.
\newblock Bayesian approach to a nonlinear inverse problem for a time-space
  fractional diffusion equation.
\newblock {\em Inverse Problems}, 34(12):125002, 2018.

\bibitem{LYY2015ANM}
J.J. Liu, M.~Yamamoto, and L.~Yan.
\newblock On the uniqueness and reconstruction for an inverse problem of the
  fractional diffusion process.
\newblock {\em Applied Numerical Mathematics}, 87:1--19, 2015.

\bibitem{YY2014CMA}
L.~Yan and F.~L. Yang.
\newblock Efficient kansa-type mfs algorithm for time-fractional inverse
  diffusion problems.
\newblock {\em Computers \& Mathematics with Applications}, 67(8):1507--1520,
  2014.

\bibitem{Iglesias+Law+Kody2013ensemble}
M.~A. Iglesias, K.~JH Law, and A.~M. Stuart.
\newblock Ensemble kalman methods for inverse problems.
\newblock {\em Inverse Problems}, 29(4):045001, 2013.

\bibitem{Iglesias2016regularizing}
M.~A. Iglesias.
\newblock A regularizing iterative ensemble kalman method for pde-constrained
  inverse problems.
\newblock {\em Inverse Problems}, 32(2):025002, 2016.

\bibitem{Law+Stuart2015data}
K.~Law, A.~Stuart, and K.~Zygalakis.
\newblock {\em Data Assimilation: A Mathematical Introduction}, volume~62.
\newblock Springer, 2015.

\bibitem{SS2017SJNA}
C.~Schillings and A.~M. Stuart.
\newblock Analysis of the ensemble kalman filter for inverse problems.
\newblock {\em SIAM Journal on Numerical Analysis}, 55(3):1264--1290, 2017.

\bibitem{SS2018AA}
C.~Schillings and A.~M. Stuart.
\newblock Convergence analysis of ensemble kalman inversion: the linear, noisy
  case.
\newblock {\em Applicable Analysis}, 97(1):107--123, 2018.

\bibitem{Evensen1994}
G.~Evensen.
\newblock Sequential data assimilation with a nonlinear quasi-geostrophic model
  using monte carlo methods to forecast error statistics.
\newblock {\em Journal of Geophysical Research: Oceans}, 99(C5):10143--10162,
  1994.

\bibitem{Asher2015review}
M.J. Asher, B.F.W. Croke, and L.J.M. Jakeman, A.J .and~Peeters.
\newblock A review of surrogate models and their application to groundwater
  modeling.
\newblock {\em Water Resources Research}, 51(8):5957--5973, 2015.

\bibitem{Li+Lin+Zhang2014adaptive}
W.~Li, G.~Lin, and D.~Zhang.
\newblock An adaptive anova-based pckf for high-dimensional nonlinear inverse
  modeling.
\newblock {\em Journal of Computational Physics}, 258:752--772, 2014.

\bibitem{Li+Zhang+Lin2015surrogate}
W.~Li, D.~Zhang, and G.~Lin.
\newblock A surrogate-based adaptive sampling approach for history matching and
  uncertainty quantification.
\newblock In {\em SPE Reservoir Simulation Symposium}. Society of Petroleum
  Engineers, 2015.

\bibitem{ju2018adaptive}
L.~Ju, J.~Zhang, L.~Meng, L.~Wu, and L.~Zeng.
\newblock An adaptive gaussian process-based iterative ensemble smoother for
  data assimilation.
\newblock {\em Advances in Water Resources}, 115:125--135, 2018.

\bibitem{YZ2018IJUQ}
L.~Yan and T.~Zhou.
\newblock An adaptive multi-fidelity pc-based ensemble kalman inversion for
  inverse problems.
\newblock {\em arXiv preprint arXiv:1809.08931}, 2018.

\bibitem{QMN2015book}
A.~Quarteroni, A.~Manzoni, and F.~Negri.
\newblock {\em Reduced basis methods for partial differential equations: an
  introduction}, volume~92.
\newblock Springer, 2015.

\bibitem{RHP2007reduced}
G.~Rozza, D.~B.~P. Huynh, and A.~T. Patera.
\newblock Reduced basis approximation and a posteriori error estimation for
  affinely parametrized elliptic coercive partial differential equations.
\newblock {\em Archives of Computational Methods in Engineering}, 15(3):1,
  2007.

\bibitem{HU2018JCP}
J.~S. Hesthaven and S.~Ubbiali.
\newblock Non-intrusive reduced order modeling of nonlinear problems using
  neural networks.
\newblock {\em Journal of Computational Physics}, 363:55--78, 2018.

\bibitem{GH2018CMAME}
M.~Guo and J.~S. Hesthaven.
\newblock Reduced order modeling for nonlinear structural analysis using
  gaussian process regression.
\newblock {\em Computer Methods in Applied Mechanics and Engineering},
  341:807--826, 2018.

\bibitem{ADN2013NMPDE}
C.~Audouze, F.~De~Vuyst, and P.~B. Nair.
\newblock Nonintrusive reduced-order modeling of parametrized time-dependent
  partial differential equations.
\newblock {\em Numerical Methods for Partial Differential Equations},
  29(5):1587--1628, 2013.

\bibitem{XFPN2017CMAME}
D.~Xiao, F.~Fang, C.~Pain, and I.~Navon.
\newblock A parameterized non-intrusive reduced order model and error analysis
  for general time-dependent nonlinear partial differential equations and its
  applications.
\newblock {\em Computer Methods in Applied Mechanics and Engineering},
  317:868--889, 2017.

\bibitem{GH2019CMAME}
M.~Guo and J.~S. Hesthaven.
\newblock Data-driven reduced order modeling for time-dependent problems.
\newblock {\em Computer Methods in Applied Mechanics and Engineering},
  345:75--99, 2019.

\bibitem{BGW2015SIRV}
P.r Benner, S.~Gugercin, and K.~Willcox.
\newblock A survey of projection-based model reduction methods for parametric
  dynamical systems.
\newblock {\em SIAM review}, 57(4):483--531, 2015.

\bibitem{BJWW2000IP}
H.T. Banks, M.~L. Joyner, B.~Wincheski, and W.~P. Winfree.
\newblock Nondestructive evaluation using a reduced-order computational
  methodology.
\newblock {\em Inverse Problems}, 16(4):929, 2000.

\bibitem{YYL2018JCP}
F.~L. Yang, L.~Yan, and L.~Ling.
\newblock Doubly stochastic radial basis function methods.
\newblock {\em Journal of Computational Physics}, 363:87--97, 2018.

\bibitem{YZ2019JCP}
L.~Yan and T.~Zhou.
\newblock Adaptive multi-fidelity polynomial chaos approach to bayesian
  inference in inverse problems.
\newblock {\em Journal of Computational Physics}, 381:110--128, 2019.

\bibitem{yan+guo2015}
L.~Yan and L.~Guo.
\newblock Stochastic collocation algorithms using $l_1$-minimization for
  {B}ayesian solution of inverse problems.
\newblock {\em SIAM Journal on Scientific Computing}, 37(3):A1410--A1435, 2015.

\bibitem{Cui2014data}
T.~Cui, Y.~M. Marzouk, and K.~Willcox.
\newblock Data-driven model reduction for the {B}ayesian solution of inverse
  problems.
\newblock {\em International Journal for Numerical Methods in Engineering},
  102(5):966--990, 2015.

\bibitem{Yan+Zhang2017IP}
L.~Yan and Y.X. Zhang.
\newblock Convergence analysis of surrogate-based methods for bayesian inverse
  problems.
\newblock {\em Inverse Problems}, 33(12):125001, 2017.

\end{thebibliography}

\end{document}